\documentclass{amsart}
\usepackage{amssymb}

\usepackage{amscd}
\usepackage{thmdefs}
\usepackage{amsmath}
\usepackage{graphicx}

\providecommand{\U}[1]{\protect\rule{.1in}{.1in}}
\theoremstyle{definition}
\theoremstyle{remark}
\numberwithin{equation}{section}

\input tcilatex

\begin{document}
\title[Graph Measures]{Measures on Graphs and Certain Groupoid Measures}
\author{Ilwoo Cho}
\address{Saint Ambrose University, Dep. of Math, 115 McM Hall, 518 W. Locust St.,
Davenport, Iowa, U. S. A.}
\email{chowoo@sau.edu}
\thanks{I really appreciate all supports from Saint Ambrose University. I also thank
Prof. T. Anderson and Prof. V. Vega for the valuable discussion.}
\thanks{}
\date{June, 2006}
\subjclass{}
\keywords{Finite Directed Graphs, Energy Measures, Diagram Measures, Graph Groupoid
Measures, Reduced Diagram Measures, Graph Integrals, Graph Hilbert Spaces,
Graph von Neumann Algebras.}
\dedicatory{}
\maketitle

\begin{abstract}
The main purpose of this paper is to introduce several measures determined
by a given finite directed graph. To construct $\sigma $-algebras for those
measures, we consider several algebraic structures induced by $G$; (i) the
free semigroupoid $\Bbb{F}^{+}(G^{\symbol{94}})$ of the shadowed graph $G^{%
\symbol{94}}$ $=$ $G$ $\cup $ $G^{-1}$ (ii) the graph groupoid $\Bbb{G}$ of $%
G$, (iii) the disgram set $D(G^{\symbol{94}})$ and (iv) the reduced diagram
set $D_{r}(G^{\symbol{94}}).$ The graph measures $\mu _{G^{\symbol{94}}}$
determined by (i) is the energy measure measuing how much energy we spent
when we have some movements on $G.$ The graph measures $\mu _{\delta }$
determined by (iii) is the diagram measure measuring how long we moved
consequently from the starting positions (which are vertices) of some
movements on $G.$ The graph measures $\mu _{\Bbb{G}}$ and $\mu _{\delta
^{r}} $ determined by (ii) and (iv) are the (graph) groupoid measure and the
(quotient-)groupoid measure, respectively.\ We show that above graph
measurings are invariants on finite directed graphs, when we consider
shadowed graphs are certain two-colored graphs. Also, we will consider the
reduced diagram measure theory on graphs. In the final chapter, we will show
that if two finite directed graphs $G_{1}$ and $G_{2}$ are graph-isomorphic,
then the von Neumann algebras $L^{\infty }(\mu _{1})$ and $L^{\infty }(\mu
_{2})$ are $*$-isomorphic, where $\mu _{1}$ and $\mu _{2}$ are the same kind
of our graph measures of $G_{1}$ and $G_{2},$ respectively.
\end{abstract}

\strut

The main purpose of this paper is to introduce certain measures induced by a
finite directed graph and to introduce certain measures induced by a
groupoid, in particular, a groupoid generated by a finite directed graph.
Measure Theory of such measures is fundamental and easily understood but it
is interesting because of that the theory is depending on combinatorial
objects (i.e., graphs).\ In fact, the operator algebra depending on such
measures looks much interesting than the measure theory. In the final
chapter of this paper, we briefly consider a von Neumann algebra $L^{\infty
} $ $(\mu _{G}),$ where $\mu _{G}$ is either one of our measures induced by
a finite directed graph or one of our graph groupoid measures. This paper
would be the first step of such combinatorial measure theoretic operator
algebra. In this paper, we will concentrate on constructing such measures
and considering their properties. We will show that our graph measuring is
quiet stable object since it is an invariant on finite directed graphs under
certain additional assumption. From this invariance, we can see that if two
finite directed graphs $G_{1}$ and $G_{2}$ are graph-isomorphic, then the
corresponding our graph measures $\mu _{G_{1}}$ and $\mu _{G_{2}}$ are
equivalent and hence the von Neumann algebras $L^{\infty }$ $(\mu _{G_{1}})$
and $L^{\infty }$ $(\mu _{G_{2}})$ are $*$-isomorphic.

\strut

Let $G$ be a finite directed graph, with its vertex set $V(G)$ and its edge
set $E(G).$ Throughout this paper, we say that a graph is finite if $\left|
V(G)\right| $ $<$ $\infty $ and $\left| E(G)\right| $ $<$ $\infty .$ Let $v$ 
$\in $ $V(G)$ be a vertex. On the graph $G,$ we can have a finite path $w$ $%
= $ $e_{1}$ $e_{2}$ ... $e_{k},$ where $e_{1},$ ..., $e_{k}$ are directed
edges in $E(G)$. In this case, we say that the edges $e_{1},$ ..., $e_{k}$
are admissible. Define the length $\left| w\right| $ of $w$ by $k$ which is
nothing but the cardinality of the admissible edges constructing the finite
path. Suppose a finite path $w$ has its initial vertex (or its source) $%
v_{1} $ and its terminal vertex (or its range) $v_{2}.$ Then we write $w$ $=$
$v_{1}$ $w$ $v_{2}$ to emphasize the initial vertex and the terminal vertex
of $w.$ Sometimes, we will denote $w$ $=$ $v_{1}$ $w$ (or $w$ $=$ $w$ $v_{2}$%
), to emphasize the initial vertex of $w$ (resp. the terminal vertex of $w$%
). If $w$ $=$ $v_{1}$ $w$ $v_{2}$ with $v_{1},$ $v_{2}$ $\in $ $V(G),$ then
we also say that $v_{1}$ and $w$ (resp. $w$ and $v_{2}$) are admissible.
Notice that every finite path $w$ can be understood as a word in $E(G),$
under the admissibility. Denote the set of all finite paths including the
empty word $\emptyset $ by $FP(G)$.

\strut

Let $w_{1}$ and $w_{2}$ be in $V(G)$ $\cup $ $FP(G)$ and assume that the
product of two words $w_{1}$ $w_{2}$ of $w_{1}$ and $w_{2}$ is again
contained in $V(G)$ $\cup $ $FP(G)$. Then we say that $w_{1}$ and $w_{2}$
are admissible. Note that, in general, even though $w_{1}$ and $w_{2}$ are
admissible, $w_{2}$ and $w_{1}$ are not admissible, since the admissibility
is totally depending on the direction on the graph $G.$

\strut

Define the free semigroupoid $\mathbb{F}^{+}(G)$ $=$ $V(G)$ $\cup $ $FP(G),$
with its binary operation ($\cdot $) called the admissibility, where $%
\emptyset $ is the empty word in $V(G)$ $\cup $ $FP(G).$ The admissibility
of on $\Bbb{F}^{+}(G)$ defined by

\strut

\begin{center}
$(w_{1},w_{2})\longmapsto\left\{ 
\begin{array}{ll}
w_{1}w_{2} & \text{if }w_{1}w_{2}\in V(G)\cup FP(G) \\ 
\emptyset & \text{otherwise.}
\end{array}
\right. $
\end{center}

\strut\strut\strut

One of the main purpose of this paper is to introduce several measures on a
finite directed graph. First, we will observe two kinds of graph measures $%
\mu _{G^{\symbol{94}}}$ and $\mu _{\delta }$ on a graph $G,$ so-called the
energy measure and the diagram measure of the shadowed graph $G^{\symbol{94}%
} $ of $G$, respectively. The construction of $\mu _{G^{\symbol{94}}}$ and $%
\mu _{\delta }$ is basically same, but their measurable spaces are
different. The energy measure $\mu _{G^{\symbol{94}}}$ measures how much
energy we spent when we moves on the shadowed graph $G^{\symbol{94}}$. And
the diagram measure $\mu _{\delta }$ measures how long we moved consequently
from the starting point (or a initial vertex of the movement) on the graph $%
G^{\symbol{94}}.$ Also, based on these measures, we consider the groupoid
measures $\mu _{\Bbb{G}}$ and $\mu _{\delta ^{r}}$ called the graph groupoid
measure and reduced diagram measure of $G^{\symbol{94}},$ respectively..

\strut

To construct the measurable spaces of the energy measure $\mu _{G^{\symbol{94%
}}},$ we will consider the free semigroupoid $\Bbb{F}^{+}(G^{\symbol{94}})$
of the shadowed graph $G^{\symbol{94}}$ $=$ $G$ $\cup $ $G^{-1}$ of the
given finite directed graph $G$, with its vertex set $V(G^{\symbol{94}})$ $=$
$V(G)$ $=$ $V(G^{-1})$ and its edge set $E(G)$ $=$ $E(G)$ $\cup $ $%
E(G^{-1}), $ where $G^{-1}$ is the shadow of $G$ that is the opposite
directed graph of $G.$ Since $G^{\symbol{94}}$ is a new directed graph, we
can construct the free semigroupoid $\Bbb{F}^{+}(G^{\symbol{94}})$ of the
shadowed graph $G^{\symbol{94}}.$ The pair $U_{1}$ $=$ $(\Bbb{F}^{+}(G^{%
\symbol{94}}),$ $P\left( \Bbb{F}^{+}(G^{\symbol{94}})\right) )$ is a
measurable space, where $P(X)$ means the power set of an arbitrary set $X.$
On this measurable space $U_{1}$, we can define a measure $\mu _{G^{\symbol{%
94}}}$ by $\mu _{G^{\symbol{94}}}$ $\overset{def}{=}$ $d$ $\cup $ $\rho $
such that

\strut

\begin{center}
$\mu _{G^{\symbol{94}}}(S)=d\left( S\text{ }\cap \text{ }V(G^{\symbol{94}%
})\right) +\rho \left( S\text{ }\cap \text{ }FP(G^{\symbol{94}})\right) ,$
\end{center}

\strut

for all $S$ $\in $ $P\left( \Bbb{F}^{+}(G^{\symbol{94}})\right) $, where $d$
is the degree measure on $V(G^{\symbol{94}})$ defined by

\strut

\begin{center}
$d(V)=0,$ \ \ for all \ $V$ $\subseteq $ $V(G^{\symbol{94}})$
\end{center}

\strut

and where $\rho $ is the length measure on $FP(G^{\symbol{94}})$ defined by

\strut

\begin{center}
$\rho (F)=\underset{w\in F}{\sum }\left| w\right| ,$ \ for all \ $F$ $%
\subseteq $ $FP(G^{\symbol{94}}).$
\end{center}

\strut

This measure $\mu _{G^{\symbol{94}}}$ is the energy measure of $G.$

\strut

To define the diagram-length measure $\mu _{\delta }$ of $G,$ we construct
the algebraic structure $D(G^{\symbol{94}})$ called the diagram set of $G^{%
\symbol{94}}.$\ Define the diagram map $\delta $ $:$ $\mathbb{F}^{+}(G^{%
\symbol{94}})$ $\rightarrow $ $\mathbb{F}^{+}(G^{\symbol{94}})$ by $w$ $%
\mapsto $ $\delta _{w},$ where the diagram $\delta _{w}$ of $w$ is nothing
but the graphical image of $w$ in $\mathbb{R}^{2}$ on $G.$ The image $D(G^{%
\symbol{94}})$ $=$ $\delta \left( \Bbb{F}^{+}(G^{\symbol{94}})\right) $ of $%
\delta $ is a subset of $\Bbb{F}^{+}(G^{\symbol{94}})$ and it also has its
admissibility, as binary operation on it, inherited by that of $\mathbb{F}%
^{+}(G^{\symbol{94}})$. i.e., $\delta _{w}$ $\delta _{w^{\prime }}$ $=$ $%
\delta _{ww^{\prime }}.$ Thus the diagram set $D(G^{\symbol{94}})$ with the
inherited admissibility is again an algebraic structure induced by $G$%
.\strut The measurable space $U_{2}$ $=$ $(D(G^{\symbol{94}}),$ $P\left(
D(G^{\symbol{94}})\right) )$ is well-determined and the diagram measure $\mu
_{\delta }$ $=$ $d$ $\cup $ $\rho $ is also well-defined by

\strut

\begin{center}
$\mu _{\delta }(S)=d\left( S\text{ }\cap \text{ }D_{V}(G^{\symbol{94}%
})\right) +\rho \left( S\text{ }\cap \text{ }D_{FP}(G^{\symbol{94}})\right)
, $
\end{center}

\strut

for all $S$ $\in $ $P\left( D(G^{\symbol{94}})\right) ,$ where $D_{V}(G^{%
\symbol{94}})$ $\overset{def}{=}$ $D(G^{\symbol{94}})$ $\cap $ $V(G^{\symbol{%
94}})$ $=$ $V(G^{\symbol{94}})$ and $D_{FP}(G^{\symbol{94}})$ $\overset{def}{%
=}$ $D(G^{\symbol{94}})$ $\cap $ $FP(G^{\symbol{94}}).$ Here, $d$ is the
vertex measure and $\rho $ is the diagram-length measure defined by $\rho $ $%
\overset{def}{=}$ $\rho \mid _{D_{FP}(G^{\symbol{94}})}.$

\strut

Based on $\mu _{G^{\symbol{94}}}$ and $\mu _{\delta },$ we will think about
certain groupoid measures. By defining the reducing relation (RR) on the
free semigroupoid $\Bbb{F}^{+}(G^{\symbol{94}}),$ we can construct the graph
groupoid $\Bbb{G}$ of $G,$ where

\strut

(RR) \ \ \ $\ \ \ \ \ \ \ \ \ \ \ \ ww^{-1}=v_{1}$ \ \ and \ $w^{-1}w=v_{2},$

\strut

for all $w$ $=$ $v_{1}$ $w$ $v_{2}$ in $\Bbb{F}^{+}(G^{\symbol{94}}),$ with $%
v_{1}$, $v_{2}$ $\in $ $V(G^{\symbol{94}}).$ (If $w$ is a vertex in $V(G^{%
\symbol{94}}),$ we can regard it as $w$ $=$ $w$ $w$ $w.$) The graph groupoid 
$\Bbb{G}$ with the same admissibility on $\Bbb{F}^{+}(G^{\symbol{94}})$ is
indeed a categroail groupoid with its objects $V(G^{\symbol{94}})$ and its
morphisms $FP_{r}(G^{\symbol{94}})$ $\overset{def}{=}$ $\Bbb{G}$ $\setminus $
$V(G^{\symbol{94}}).$ The graph groupoid measure $\mu _{\Bbb{G}}$ is nothing
but $\mu _{G^{\symbol{94}}}\mid _{\Bbb{G}}$ with respect to the measurable
space $U_{3}$ $=$ $(\Bbb{G},$ $P\left( \Bbb{G}\right) ).$

\strut

Finally, define a reduced diagram set $D_{r}(G^{\symbol{94}})$ by the image
of the reduced diagram map $\delta ^{r}$ $:$ $\Bbb{F}^{+}(G^{\symbol{94}})$ $%
\rightarrow $ $\Bbb{F}^{+}(G^{\symbol{94}})$ by $w$ $\mapsto $ $\delta
_{w}^{r},$ where $\delta _{w}^{r}$ is the reduced diagram of $w.$ i.e., an
element $\delta _{w}^{r}$ of $w$ is a diagram $\delta _{w}$ under (RR).
i.e., $\delta ^{r}$ $=$ $\delta ^{r}$ $\circ $ $\delta .$ The image $%
D_{r}(G^{\symbol{94}})$ $\overset{def}{=}$ $\delta ^{r}\left( \Bbb{F}^{+}(G^{%
\symbol{94}})\right) $ $=$ $\delta ^{r}\left( D(G^{\symbol{94}})\right) $ is
called a reduced diagram set and it has its inherited admissibility from $%
\Bbb{G}.$ It is easy to check that the reduced diagram set $D_{r}(G^{\symbol{%
94}})$ is a sub-structure of the graph groupoid $\Bbb{G}$ of $G.$ We can
define the reduced diagram measure $\mu _{\delta ^{r}}\Bbb{\ }$by $\mu _{%
\Bbb{G}}\mid _{D_{r}(G^{\symbol{94}})}$ with respect to the measurable space 
$U_{4}$ $=$ $(D_{r}(G^{\symbol{94}}),$ $P\left( D_{r}(G^{\symbol{94}%
})\right) ).$

\strut

As a main result of this paper, we will show that such graph measurings are
invariants on finite directed graphs under an additional assumption. i.e.,
the finite directed graphs $G_{1}$ and $G_{2}$ are graph-isomorphic if and
only if the corresponding shadowed graphs $G_{1}^{\symbol{94}}$ and $G_{2}^{%
\symbol{94}}$ are two-colored-graph-isomorphic if and only if the
corresponding graph measures are equivalent. To study algebraic structures
(free semigroupoids, diagram sets, graph groupoids, and reduced diagram
sets), itself would be also very interesting, since they are combinatorial
object depending group-like structures which are not groups. For instance,
we can understand diagram sets (resp. reduced diagram sets) as quotient
structures of free semigroupoids (resp. graph groupoids).

\strut

With respect to those measures on graphs, we will consider Measure Theory.
In particular, we will concentrate on observing Integration on the reduced
diagram measure space $(D_{r}(G^{\symbol{94}}),$ $P\left( D_{r}(G^{\symbol{94%
}})\right) ,$ $\mu _{\delta ^{r}}).$

\strut

In the final chapter, we introduce von Neumann algebras induced by our graph
measures. We can see that the graph von Neumann algebra $L^{\infty }(\mu
_{G_{1}})$ and $L^{\infty }(\mu _{G_{2}})$ are $*$-isomorphic if the finite
directed graphs $G_{1}$ and $G_{2}$ are graph-isomorphic, where $\mu
_{G_{k}} $ are the same kind of our graph measures of $G_{k},$ for $k$ $=$ $%
1,$ $2.$

\strut \strut

\strut

\strut \strut

\section{Measures Induced by a Graph}

\strut

\strut

In this chapter, we will construct measures induced by a finite directed
graph and observe their properties. For convenience, we will call such
measures graph measures. We will define four kinds of graph measures. Two of
them are measures measuring certain activity on a graph (the energy measure
and the diagram measure)\ and the other two measures contain the measure
theoretic information of a groupoid generated by the graph (the graph
groupoid measure and the reduced diagram measure).

\strut

Throughout this chapter, let $G$ be a finite directed graph with its vertex
set $V(G)$ and its edge set $E(G)$. We say that edges $e_{1},$ ..., $e_{k}$
are admissible (or connected in the order $(e_{1},$ ..., $e_{k})$) if there
exists a finite path $w$ such that $w$ $=$ $e_{1}$ ... $e_{k}.$ In this
case, the length $\left| w\right| $ of this finite path $w$ is defined to be 
$k$ which is the cardinality of directed edges constructing the finite path $%
w.$ Denote the set of all finite paths by $FP(G).$ Clearly, the edge set $%
E(G)$ of $G$ is contained in $FP(G),$ since all edges are regarded as finite
paths with their lengths $1.$ Suppose $w$ is an arbitrary finite path in $%
FP(G)$ and assume that it has its initial vertex (or a source) $v_{1}$ and
its terminal vertex (or a range) $v_{2}.$ Then we denote $w$ by $w$ $=$ $%
v_{1}$ $w$ $v_{2}$ to emphasize its initial and terminal vertices. If a
finite path $w$ satisfies $w$ $=$ $v_{1}$ $w$ $v_{2},$ then we also say that 
$v_{1}$ and $w$ are admissible and that $w$ and $v_{2}$ are admissible. Let $%
w_{1}$ and $w_{2}$ be finite paths in $FP(G)$ and suppose that there exists
a finite path $w$ $=$ $w_{1}$ $w_{2}$ in $FP(G).$ Then the elements $w_{1}$
and $w_{2}$ are said to be admissible. Remark that even though $w_{1}$ and $%
w_{2}$ are admissible, in general, $w_{2}$ and $w_{1}$ are not admissible.
For instance if $w_{1}$ $=$ $v_{1}$ $w_{1}$ $v_{2}$ and $w_{2}$ $=$ $v_{2}$ $%
w_{2}$ $v_{3}$ with $v_{1}$, $v_{2},$ $v_{3}$ $\in $ $V(G)$ such that $v_{1}$
$\neq $ $v_{3},$ then $w_{1}$ and $w_{2}$ are admissible but $w_{2}$ and $%
w_{1}$ are not admissible. Notice that all finite paths in $FP(G)$ can be
regarded as words in $E(G)$ under the admissibility. So, we naturally assume
that the empty word $\emptyset $ is contained in $FP(G).$

\strut

\strut

\strut

\subsection{Free Semigroupoids and Diagram Sets of Shadowed Graphs}

\strut

\strut

Let $G$ be a finite directed graph with its vertex set $V(G)$ and its edge
set $E(G).$ Define a shadow $G^{-1}$ of $G$ by the opposite directed graph
of $G$ with its vertex set $V(G^{-1})$ $=$ $V(G)$ and its edge set $%
E(G^{-1}) $ $=$ $\{e^{-1}$ $:$ $e$ $\in $ $E(G)\}.$ Here the symbol $e^{-1}$
of an edge $e$ of $G$ means the opposite directed edge of $e.$ i.e., if $e$ $%
=$ $v$ $e$ $v^{\prime }$ with $v,$ $v^{\prime }$ $\in $ $V(G),$ the shadow$\
e^{-1}$ of $e$ satisfies $e^{-1}$ $=$ $v_{2}$ $e^{-1}$ $v_{1}.$

\strut

\begin{definition}
Let $G$ be a finite directed graph and $G^{-1},$ the shadow of $G.$ Define
the shadowed graph $G^{\symbol{94}}$ of $G$ by a finite directed graph $G^{%
\symbol{94}}$ $=$ $G$ $\cup $ $G^{-1}$ with

\strut 

$\ \ \ \ \ \ \ \ \ \ \ \ \ \ \ \ \ \ \ \ \ \ V(G^{\symbol{94}})$ $=$ $V(G)$ $%
=$ $V(G^{-1})$

and

$\ \ \ \ \ \ \ \ \ \ \ \ \ \ \ \ \ \ \ \ \ \ E(G^{\symbol{94}})=E(G)$ $\cup $
$E(G^{-1}).$
\end{definition}

\strut

For the new finite directed graph $G^{\symbol{94}}$ induced by $G,$ we can
consider its free semigroupoid $\Bbb{F}^{+}(G^{\symbol{94}}).$ Notice that
the finite path set $FP(G^{\symbol{94}})$ of the shadowed graph $G^{\symbol{%
94}}$ of $G$ does not satisfy $FP(G^{\symbol{94}})$ $=$ $FP(G)$ $\cup $ $%
FP(G^{-1}).$ In fact,

\strut

\begin{center}
$FP(G^{\symbol{94}})\supsetneqq $ $FP(G)$ $\cup $ $FP(G^{-1}),$
\end{center}

\strut

whenever $\left| E(G^{\symbol{94}})\right| $ $\geq $ $1.$ For example, if $e$
$\in $ $E(G),$ then there exists a finite path $ee^{-1}$ in $FP(G^{\symbol{94%
}}),$ but this finite path $ee^{-1}$ is neither in $FP(G)$ nor $FP(G^{-1}).$

\strut

Now, define a map $\delta $ from the free semigroupoid $\Bbb{F}^{+}(G^{%
\symbol{94}})$ of $G^{\symbol{94}}$ into itself by

\strut

\begin{center}
$\delta $ $:$ $\Bbb{F}^{+}(G^{\symbol{94}})$ $\rightarrow $ $\Bbb{F}^{+}(G^{%
\symbol{94}}),$ \ \ $w$ $\longmapsto $ $\delta _{w},$ $\forall $ $w$ $\in $ $%
\Bbb{F}^{+}(G^{\symbol{94}}),$
\end{center}

\strut

where $\delta _{w}$ means the graphical image of $w$ in $\Bbb{R}^{2}$ on $G.$
The map $\delta $ is said to be the diagram map and the image $\delta _{w}$
of $w$ is called the diagram of $w,$ for all $w$ $\in $\thinspace $\Bbb{F}%
^{+}(G^{\symbol{94}}).$

\strut

\begin{definition}
Let $G$ be a finite directed graph and $G^{\symbol{94}},$ the shadowed graph
of $G$ and let $\Bbb{F}^{+}(G^{\symbol{94}})$ be the free semigroupoid of $%
G^{\symbol{94}}.$ Suppose $\delta $ is the diagram map on $\Bbb{F}^{+}(G^{%
\symbol{94}}).$ The image $\delta \left( \Bbb{F}^{+}(G^{\symbol{94}})\right) 
$ of $\delta $ in $\Bbb{F}^{+}(G^{\symbol{94}})$ is called the diagram set
and we denote it by $D(G^{\symbol{94}}).$ Let $w$ be an element in $\Bbb{F}%
^{+}(G^{\symbol{94}})$ satisfying that $\delta _{w}$ $=$ $w.$ Then $w$ is
said to be basic in $\Bbb{F}^{+}(G^{\symbol{94}}).$ Under the inherited
admissibility from $\Bbb{F}^{+}(G^{\symbol{94}}),$ the diagram set $D(G^{%
\symbol{94}})$ is a substructure of $\Bbb{F}^{+}(G^{\symbol{94}}).$ i.e., $%
\delta _{w_{1}}$ $\cdot $ $\delta _{w_{2}}$ $=$ $\delta _{w_{1}w_{2}},$ for
all $\delta _{w_{1}},$ $\delta _{w_{2}}$ $\in $ $D(G^{\symbol{94}}).$ We
will denote the diagram set $(D(G^{\symbol{94}}),$ $\cdot )$ simply by $D(G^{%
\symbol{94}}).$
\end{definition}

\strut

By definition, we can re-define the set $D(G^{\symbol{94}})$ by the subset
of $\Bbb{F}^{+}(G^{\symbol{94}})$ consisting of all basic elements. It is
easy to see that all vertices and all edges are basic elements. Suppose $l$ $%
=$ $v$ $l$ $v$ is a loop finite path in $FP(G^{\symbol{94}}),$ with $v$ $\in 
$ $V(G^{\symbol{94}}),$ and assume that its diagram $\delta _{l}$ is
identified with $w_{0}$ in $\Bbb{F}^{+}(G^{\symbol{94}}).$ Then there exists 
$n$ $\in $ $\Bbb{N}$ such that $l$ $=$ $w_{0}^{n}.$ If $n$ $=$ $1,$ then $l$ 
$=$ $w_{0}$ is the basic element in $\Bbb{F}^{+}(G^{\symbol{94}}).$ In fact,
for any $n$ $\in $ $\Bbb{N},$ the diagrams $\delta _{l^{n}}$ $=$ $w_{0},$ in 
$D(G^{\symbol{94}}).$

\strut

Consider the admissibility on $D(G^{\symbol{94}}).$ If $\delta _{1}$ and $%
\delta _{2}$ are diagrams in $D(G^{\symbol{94}}),$ then there exist basic
elements $w_{1}$ and $w_{2}$ such that $w_{1}$ $=$ $\delta _{1}$ and $w_{2}$ 
$=$ $\delta _{2},$ respectively. Suppose that $w_{1}$ $w_{2}$ $\neq $ $%
\emptyset .$ Then we have that $\delta _{1}$ $\delta _{2}$ $=$ $\delta
_{w_{1}w_{2}}.$ Otherwise, $\delta _{1}$ $\delta _{2}$ $=$ $\delta
_{\emptyset }$ $=$ $\emptyset .$

\strut \strut \strut \strut

\begin{proposition}
Define a relation $\mathcal{R}$ on $\Bbb{F}^{+}(G^{\symbol{94}})$ by

\strut 

$\ \ \ \ \ \ \ \ \ w_{1}$ $\mathcal{R}$ $w_{2}$ $\overset{def}{%
\Longleftrightarrow }$ $\delta _{w_{1}}=\delta _{w_{2}},$ \ for \ $w_{1},$ $%
w_{2}$ $\in $ $\Bbb{F}^{+}(G^{\symbol{94}}).$

\strut 

Then the relation $\mathcal{R}$ is an equivalent relation on $\Bbb{F}^{+}(G^{%
\symbol{94}})$. $\square $
\end{proposition}

\strut

So, the substructure $D(G^{\symbol{94}})$ of $\Bbb{F}^{+}(G^{\symbol{94}})$
can be understood as a quotient algebraic structure $\Bbb{F}^{+}(G^{\symbol{%
94}})$ $/$ $\mathcal{R}$, with its quotient map $\delta .$ Also, diagrams
are understood as equivalence classes.

\strut \strut \strut \strut

\strut

\subsection{Graph Groupoids and Reduced Diagram Sets}

\strut

\strut

In this section, we will define graph groupoids and reduced diagram sets, by
defining so-called the reducing relation. The purpose of this section is to
construct certain groupoid measures induced by graphs.

\strut

\begin{definition}
Let $G$ be a finite directed graph and $G^{\symbol{94}},$ the shadowed graph
of $G$ and let $\Bbb{F}^{+}(G^{\symbol{94}})$ be the free semigroupoid of $%
G^{\symbol{94}}.$ Define the reducing relation (RR) on $\Bbb{F}^{+}(G^{%
\symbol{94}})$ by

\strut 

(RR) \ \ \ \ \ \ \ \ \ \ \ \ $ww^{-1}=v_{1}$ and $w^{-1}w=v_{2},$

\strut 

for all $w$ $=$ $v_{1}$ $w$ $v_{2}$ in $\Bbb{F}^{+}(G^{\symbol{94}}),$ with $%
v_{1},$ $v_{2}$ $\in $ $V(G^{\symbol{94}}).$ (We can regard a vertex $v$ by $%
v$ $v$ $v.$) The graph groupoid $\Bbb{G}$ of $G$ is defined by the set $\Bbb{%
F}_{r}^{+}(G^{\symbol{94}})$ with the same admissibility on $\Bbb{F}^{+}(G^{%
\symbol{94}})$ under (RR). i.e., $\Bbb{G}$ $=$ $(\Bbb{F}_{r}^{+}(G^{\symbol{%
94}}),$ $\cdot ).$ Finite paths in $\Bbb{G}$ are called reduced finite paths.
\end{definition}

\strut

Let $w$ $=$ $v_{1}$ $w$ $v_{2}$ be a finite path in $FP(G^{\symbol{94}}),$
with $v_{1},$ $v_{2}$ $\in $ $V(G^{\symbol{94}}).$ Then both $ww^{-1}$ and $%
w^{-1}w$ are finite paths contained in $FP(G^{\symbol{94}})$ $\subset $ $%
\Bbb{F}^{+}(G^{\symbol{94}}).$ However, these elements $ww^{-1}$ and $%
w^{-1}w $ are identified with vertices $v_{1}$ and $v_{2}$ in $\Bbb{F}%
_{r}^{+}(G^{\symbol{94}})$ under the reducing relation (RR). More generally,
if we have a finite path $w$ $=$ $w_{1}$ ... $w_{k}$ $w_{k}^{-1}$ ... $w_{n}$
$\neq $ $\emptyset $ in $FP(G^{\symbol{94}}),$ then this element $w$ is
identified with $w^{\prime }$ $=$ $w_{1}$ ... $w_{k-1}$ $w_{k+1}$ ... $w_{n}$
in $\Bbb{G}.$ Notice that since $w$ $\neq $ $\emptyset ,$ $w_{k-1}$ and $%
w_{k+1}$ are admissible via $w_{k}$ $w_{k}^{-1}$ which is a vertex in $\Bbb{G%
},$ and hence $w^{\prime }$ $\neq $ $\emptyset $ is well-determined both in $%
\Bbb{F}^{+}(G^{\symbol{94}})$ and $\Bbb{G}.$ For convenience, we denote the
reduced finite path set in $\Bbb{G}$ by $FP_{r}$ $(G^{\symbol{94}}).$ i.e., $%
FP_{r}$ $(G^{\symbol{94}})$ $\overset{def}{=}$ $\Bbb{F}_{r}^{+}$ $(G^{%
\symbol{94}})$ $\setminus $ $V(G^{\symbol{94}}).$ Notice that the graph
groupoid $\Bbb{G}$ of $G$ is indeed a categorial groupoid with its objects $%
V(G^{\symbol{94}})$ and its morphisms $FP_{r}(G^{\symbol{94}}).$

\strut

For the graph groupoid $\Bbb{G}$ of $G,$ we define a map $\delta ^{r}$ $:$ $%
\Bbb{G}$ $\rightarrow $ $\Bbb{G}$ by $\delta \mid _{\Bbb{G}},$ where $\delta 
$ is the diagram map defined in Section 1.1. i.e., the map $\delta ^{r}$ is
the diagram map on $\Bbb{G}.$ The map $\delta ^{r}$ is said to be the
reduced diagram map on $\Bbb{G}$ and the images $\delta _{w}^{r}$ of
elements $w$ $\in $ $\Bbb{G}$ are called the reduced diagrams of $w.$

\strut

\begin{definition}
Let $G$ be a finite directed graph and $\Bbb{G},$ the corresponding graph
groupoid of $G$ and let $\delta ^{r}$ be the reduced diagram map on $\Bbb{G}.
$ Then the image $D_{r}(G^{\symbol{94}})$ $\overset{def}{=}$ $\delta
^{r}\left( \Bbb{G}\right) $ is called the reduced diagram set with the
inherited admissibility on $\Bbb{G}.$ If an element $w$ $\in $ $\Bbb{G}$
satisfies $\delta _{w}^{r}$ $=$ $w,$ then we say that $w$ is a basic element
in $\Bbb{G}.$
\end{definition}

\strut

Recall that the diagram set $D(G^{\symbol{94}})$ can be regarded as a
quotient structure $\Bbb{F}^{+}(G^{\symbol{94}})$ $/$ $\mathcal{R}$ of the
free semigroupoid $\Bbb{F}^{+}(G^{\symbol{94}}),$ by considering the diagram
map $\delta $ as a quotient map providing an equivalence relation on $\Bbb{F}%
^{+}(G^{\symbol{94}}).$ Similarly, under the reducing relation (RR), we can
understand the reduced diagram set $D_{r}(G^{\symbol{94}})$ is a quotient
structure of the graph groupoid $\Bbb{G}$ with its quotient map $\delta
^{r}. $

\strut

\begin{proposition}
Define a relation $\mathcal{R}_{r}$ on the graph groupoid $\Bbb{G}$ of $G$ by

\strut 

$\ \ \ \ \ \ \ \ \ \ \ w_{1}$ $\mathcal{R}_{r}$ $w_{2}\overset{def}{%
\Longleftrightarrow }\delta _{w_{1}}^{r}=\delta _{w_{2}}^{r},$ for all $%
w_{1},$ $w_{2}$ $\in $ $\Bbb{G}.$

\strut 

Then the relation $\mathcal{R}_{r}$ is an equivalence relation on $\Bbb{G}.$ 
$\square $
\end{proposition}

\strut

\strut \strut

\strut

\subsection{Energy Measures and Diagram Measures on Graphs}

\strut

\strut

\strut

In this section, we will define two kinds of measures induced by a finite
directed graph so-called the energy measures and diagram measures. To do
that we will define certain measurable spaces for them. Throughout this
section, we also let $G$ be a finite directed graph and $G^{\symbol{94}},$
the shadowed graph of $G$ and let $\Bbb{F}^{+}(G^{\symbol{94}})$ be the free
semigroupoid of $G^{\symbol{94}}$ and $D(G^{\symbol{94}}),$ the diagram set
of $G^{\symbol{94}}.$ We define the energy measure $\mu _{G^{\symbol{94}}}$
of $G$ in the following definition. This measure $\mu _{G^{\symbol{94}}}$
measures how much energy we spent when we make some movements on the graph $%
G.$ The movement on $G$ is represented by a finite path. The energy measure $%
\mu _{G^{\symbol{94}}}$ consists of the vertex measure $d$ and the finite
path measure $\rho .$ Throughout this paper, we will let the vertex measure $%
d$ be the constant function $0.$ The finite path measure $\rho $ measures
the energy to move along the routes (finite paths). For convenience, we will
assume that all vertices have no energy and that all edges have the same
quantities of energy. However, by giving the weights to edges, we can
consider the case when edges have the different quantities. We will not
consider such weighted case in this paper.

\strut

\begin{definition}
Let $G$ be the given graph and let $\Bbb{F}^{+}(G^{\symbol{94}})$ be the
free semigroupoid of the shadowed graph $G^{\symbol{94}}$ of $G.$ The pair $(%
\Bbb{F}^{+}(G^{\symbol{94}}),$ $P\left( \Bbb{F}^{+}(G^{\symbol{94}})\right) )
$ is called the energy measurable space of $G,$ where $P(X)$ means the power
set of an arbitrary set $X.$ Define the energy measure $\mu _{G^{\symbol{94}%
}}$ of $G$ by a measure determined on the $\sigma $-algebra $P\left( \Bbb{F}%
^{+}(G^{\symbol{94}})\right) $ satisfying that $\mu _{G^{\symbol{94}}}$ $=$ $%
d$ $\cup $ $\rho $ such that

\strut 

(1.1)$\ \ \ \ \ \ \ \ \mu _{G^{\symbol{94}}}\left( S\right) =d\left( S\text{ 
}\cap \text{ }V(G^{\symbol{94}})\right) +\rho \left( S\text{ }\cap \text{ }%
FP(G^{\symbol{94}})\right) ,$

\strut 

for all $S$ $\in $ $P\left( \Bbb{F}^{+}(G^{\symbol{94}})\right) ,$ where $d$
is the vertex measure on $P\left( V(G^{\symbol{94}})\right) $ such that

\strut 

(1.2) \ \ \ \ \ \ \ \ \ \ \ \ \ \ \ \ \ \ \ $d(V)=0,$ \ for all \ $V$ $%
\subseteq $ $V(G^{\symbol{94}}),$

\strut 

and where $\rho $ is the length measure on $P\left( FP(G^{\symbol{94}%
})\right) $ such that

\strut 

(1.3) \ \ \ $\ \ \ \ \ \ \ \ \ \ \ \rho \left( F\right) =\underset{w\in F}{%
\sum }\left| w\right| ,$ \ for all \ $F$ $\subseteq $ $FP(G^{\symbol{94}}).$
\end{definition}

\strut \strut \strut

\begin{remark}
We can regard the length measure $\rho $ as a weighted length measure, by
giving weights $\lambda _{e}$ $\in $ $\Bbb{R}^{+}$ to edges $e$ $\in $ $E(G^{%
\symbol{94}}).$ In this case, a weight of a finite path $w$ $=$ $e_{1}$ ... $%
e_{k}$ with admissible edges $e_{1},$ ... $e_{k}$ can be defined by $%
\sum_{i=1}^{k}$ $\lambda _{e_{i}}.$ The unweighted case is understood as the
weighted case when all weights are $1$. Also, instead of giving the zero
measure for the vertex measure, we can put some nonzero measures for the
vertices. For instance, suppose each vertex $v$ has its nonzero measure $%
\lambda _{v}$. Then the vertex measure $d(V)$ of $V$ $\subset $ $V(G^{%
\symbol{94}})$ is $\underset{v\in V}{\sum }$ $\lambda _{v}$ $\neq $ $0,$
whenever $V$ is nonempty in $V(G^{\symbol{94}}).$ Also, in this case, if $w$ 
$=$ $e_{1}$ ... $e_{k}$ $\neq $ $\emptyset $ in $FP(G^{\symbol{94}})$ and $%
e_{j}$ $=$ $v_{j}$ $e_{j}$ $v_{j}^{\prime },$ for $j$ $=$ $1,$ ..., $k,$
then the length measure $\rho $ might be

\strut 

$\ \ \ \ \ \ \ \ \ \ \ \ \ \rho $ $(\{w\})$ $=$ $\sum_{j=1}^{k}$ $\lambda
_{e_{j}}$ $+$ $\sum_{j=1}^{k}$ $(\lambda _{v_{j}}$ $+$ $\lambda
_{v_{j}}^{\prime })$.

\strut 

Remark that the vertices $v_{j}$ and $v_{j}^{\prime }$ are not necessarily
distinct in $V(G^{\symbol{94}}).$ However, by the purpose of energy
measures, it is reasonable to put the vertex measure $d$ to the zero measure.
\end{remark}

\strut

The triple $(\Bbb{F}^{+}(G^{\symbol{94}}),$ $P\left( \Bbb{F}^{+}(G^{\symbol{%
94}})\right) ,$ $\mu _{G^{\symbol{94}}})$ is a measure space. We denote this
triple simply by $(G^{\symbol{94}},$ $\mu _{G^{\symbol{94}}})$.

\strut

Now, we will define the diagram measure $\mu _{\delta }$ on a measurable
space $(D(G^{\symbol{94}}),$ $P\left( D(G^{\symbol{94}})\right) ).$ This
measure $\mu _{\delta }$ measures how long distance we moved consequently
from the starting position (a vertex) on the graph $G$.

\strut

\begin{definition}
Let $G$ be a finite directed graph and $D(G^{\symbol{94}})$ be the diagram
set of the shadowed graph $G^{\symbol{94}}$ of $G.$ Define the diagram
measure $\mu _{\delta }$ by a measure $\mu _{\delta }$ on the $\sigma $%
-algebra $P\left( D(G^{\symbol{94}})\right) $ by $\mu _{\delta }$ $=$ $d$ $%
\cup $ $\rho ,$ where

\strut 

(1.4) $\ \ \ \ \ \ \ \ \ \ \ \ \ \ \ \ \ \ d$ $=$ $0$ \ \ and \ $\rho =\rho
\mid _{D(G^{\symbol{94}})\cap FP(G^{\symbol{94}})}.$

\strut 

i.e., the diagram measure $\mu _{\delta }$ is the restriction $\mu _{G^{%
\symbol{94}}}\mid _{P(D(G^{\symbol{94}}))}$ of the energy measure $\mu _{G^{%
\symbol{94}}}$ $.$
\end{definition}

\strut \strut \strut \strut \strut \strut

\begin{example}
Let $G_{N}$ be the one-vertex-$N$-loop-edge graph. i.e., the graph $G_{N}$
is the finite directed graph with $V(G_{N})$ $=$ $\{v\}$ and $E(G_{N})$ $=$ $%
\{l_{k}$ $=$ $v$ $l_{k}$ $v$ $:$ $k$ $=$ $1,$ ..., $N\}.$ The shadowed graph 
$G_{N}^{\symbol{94}}$ is the directed graph with $V(G_{N}^{\symbol{94}})$ $=$
$\{v\}$ and $E(G_{N}^{\symbol{94}})$ $=$ $\{l_{1}^{\pm 1},$ ..., $l_{N}^{\pm
1}\}.$ So, we can have the finite path set

\strut 

$\ \ \ \ \ FP(G_{N}^{\symbol{94}})=\cup _{m=1}^{\infty }\left\{
l_{k_{1}}^{n_{1}}...l_{k_{m}}^{n_{m}}\left| 
\begin{array}{c}
k_{i}\in \{1,...,N\},\text{ }n_{i}\in \Bbb{Z}\setminus \{0\} \\ 
i=1,\text{ }2,\text{ ..., }m,\text{ \ \ }m\in \Bbb{N}
\end{array}
\right. \right\} .$

\strut 

We can easily compute that the energy measure of $\{l_{k_{1}}^{n_{1}}$ ... $%
l_{k_{m}}^{n_{m}}\}$ $\subset $ $FP(G_{N}^{\symbol{94}})$ is

\strut 

$\ \ \ \ \ \ \ \ \ \mu _{G_{N}^{\symbol{94}}}\left(
\{l_{k_{1}}^{n_{1}}...l_{k_{m}}^{n_{m}}\}\right) =\rho \left(
\{l_{k_{1}}^{n_{1}}...l_{k_{m}}^{n_{m}}\}\right) =\sum_{j=1}^{m}\left|
n_{j}\right| .$

\strut 

Let $D_{FP}(G_{N}^{\symbol{94}})$ be the subset $D(G_{N}^{\symbol{94}})$ $%
\setminus $ $(V(G_{N}^{\symbol{94}})$ $\cup $ $\{\emptyset \}).$ Then

\strut 

$\ \ \ D_{FP}(G_{N}^{\symbol{94}})=\cup _{m=1}^{N}\left\{
l_{k_{1}}^{r_{1}}...l_{k_{m}}^{r_{m}}\left| 
\begin{array}{c}
k_{i}\in \{1,...,N\},\text{ }r_{i}\in \{\pm 1\} \\ 
i=1,...,m,\text{ }m\in \Bbb{N}
\end{array}
\right. \right\} .$

\strut 

So, the diagram measure of $\{l_{k_{1}}^{r_{1}}$ $...$ $l_{k_{m}}^{r_{m}}\}$ 
$\subset $ $D_{FP}(G_{N}^{\symbol{94}})$ is

\strut 

$\ \ \ \ \ \ \ \ \ \ \ \ \ \ \ \mu _{\delta }\left(
\{l_{k_{1}}^{r_{1}}...l_{k_{m}}^{r_{m}}\}\right) =\rho \left(
\{l_{k_{1}}^{r_{1}}...l_{k_{m}}^{r_{m}}\}\right) =m.$

\strut 

Remark that the diagram $\delta \left(
l_{k_{1}}^{n_{1}}...l_{k_{m}}^{n_{m}}\right) $ of $l_{k_{1}}^{n_{1}}$ ... $%
l_{k_{m}}^{n_{m}}$ is identical to the basic element $l_{k_{1}}^{r_{1}}$ ... 
$l_{k_{m}}^{r_{m}},$ where $r_{j}$ $=$ $1$ if $n_{j}$ $>$ $0$, and $r_{j}$ $=
$ $-1$ if $n_{j}$ $<$ $0,$ for $j$ $=$ $1,$ ..., $m.$ Notice that the free
semigroupoid $\Bbb{F}^{+}(G_{N}^{\symbol{94}})$ is an infinite set, but the
diagram set $D(G_{N}^{\symbol{94}})$ is a finite set. This shows that the $%
\sigma $-algebra $P\left( \Bbb{F}^{+}(G_{N}^{\symbol{94}})\right) $ for the
energy measure $\mu _{G_{N}^{\symbol{94}}}$ consists of infinitely many
elements and the $\sigma $-algebra $P\left( D(G_{N}^{\symbol{94}})\right) $
for the diagram measure $\mu _{\delta }$ consists of finitely many elements.
And hence the energy measure $\mu _{G_{N}^{\symbol{94}}}$ is not a bounded
measure but the diagram measure $\mu _{\delta }$ is a bounded measure.
Recall that we say that a measure $\mu $ is bounded if $\mu (S)$ $<$ $\infty
,$ for all elements $S$ in the $\sigma $-algebra for $\mu .$ In fact, in our
case, it suffices to observe that

\strut 

$\ \ \ \ \ \ \ \ \ \ \ \mu _{G_{N}^{\symbol{94}}}\left( \Bbb{F}^{+}(G_{N}^{%
\symbol{94}})\right) $ $=$ $\infty $ \ and $\mu _{\delta }\left( D(G_{N}^{%
\symbol{94}})\right) $ $<$ $\infty $.

\strut 

Even though the energy measure $\mu _{G_{N}^{\symbol{94}}}$ of $G_{N}$ is
not bounded, this measure is locally bounded in the sense that $\mu _{G_{N}^{%
\symbol{94}}}(S)$ $<$ $\infty ,$ for all finite subsets $S$ of $\Bbb{F}%
^{+}(G_{N}^{\symbol{94}}).$
\end{example}

\strut

\begin{example}
Let $C_{N}$ be the one-flow circulant graph. i.e., the graph $C_{N}$ is a
finite directed graph with $V(C_{N})$ $=$ $\{v_{1},$ ..., $v_{N}\}$ and $%
E(C_{N})$ $=$ $\{e_{j}$ $=$ $v_{j}$ $e_{j}$ $v_{j+1}$ $:$ $j$ $=$ $1,$ ..., $%
N,$ $v_{N+1}$ $\overset{def}{=}$ $v_{1}\}.$ The shadowed graph $C_{N}^{%
\symbol{94}}$ of $C_{N}$ is the finite directed graph with $V(C_{N}^{\symbol{%
94}})$ $=$ $\{v_{1},$ ..., $v_{N}\}$ and $E(C_{N}^{\symbol{94}})$ $=$ $%
\{e_{1}^{\pm 1},$ ..., $e_{N}^{\pm 1}\}.$ The finite path set of the
shadowed graph $C_{N}^{\symbol{94}}$ of $C_{N}$ is

\strut \strut 

$\underset{m_{1}\geq m_{2}\in \{0,1,2,...\}}{\cup }\left\{
e_{[i]}^{r_{1}}e_{[i+1]}^{r_{1}}...e_{[i+m_{1}]}^{r_{1}}e_{[i+m_{1}]}^{r_{2}}...e_{[i+m_{1}-m_{2}]}^{r_{2}}\left| 
\begin{array}{c}
i=1,...,N \\ 
\lbrack n]\in \Bbb{Z}_{N}\setminus \{[0]\} \\ 
r_{1}\neq r_{2}\in \{\pm 1\}
\end{array}
\right. \right\} .$

\strut 

where $\Bbb{Z}_{N}$ $\overset{def}{=}$ $\{[0],$ $[1],$ ..., $[N]\}$ with its
addition $[n_{1}]$ $+$ $[n_{2}]$ $=$ $[n_{1}$ $+$ $n_{2}],$ where $[n]$ is
the equivalence class of all numbers $k$ $\in $ $\Bbb{Z}$ satisfying that $k$
$\equiv $ $n$ $(\func{mod}$ $N).$ Then the energy measure of $%
\{e_{[i]}^{r_{1}}$ $e_{[i+1]}^{r_{1}}$ $...$ $e_{[i+m_{1}]}^{r_{1}}$ $%
e_{[i+m_{1}]}^{r_{2}}$ $...$ $e_{[i+m_{1}-m_{2}]}^{r_{2}}\}$ $\subset $ $%
FP(C_{N}^{\symbol{94}})$ is computed by

\strut 

$\ \ \ \ \ \ \ \ \ \mu _{C_{N}^{\symbol{94}}}\left(
\{e_{[i]}^{r_{1}}e_{[i+1]}^{r_{1}}...e_{[i+m_{1}]}^{r_{1}}e_{[i+m_{1}]}^{r_{2}}...e_{[i+m_{1}-m_{2}]}^{r_{2}}\}\right) 
$ $=$ $m_{1}+m_{2}.$

\strut 

The set $D_{FP}(C_{N}^{\symbol{94}})$ $\overset{def}{=}$ $D(C_{N}^{\symbol{94%
}})$ $\setminus $ $(V(C_{N}^{\symbol{94}})$ $\cup $ $\{\emptyset \})$ is
determined by

\strut 

$\ \ \ \ \ D_{FP}(C_{N}^{\symbol{94}})=\cup _{m=0}^{N-1}\left\{ \left. 
\begin{array}{c}
e_{[i]}e_{[i+1]}\text{ }...\text{ }e_{[i+m]}, \\ 
e_{[i+m]}^{-1}e_{[i+m-1]}^{-1}...e_{[i]}^{-1}
\end{array}
\right. \left| 
\begin{array}{c}
i=1,...,N \\ 
\lbrack n]\in \Bbb{Z}_{N}
\end{array}
\right. \right\} .$

\strut 

So, the diagram measure of $\{e_{[i]}^{r_{1}}$ $e_{[i+1]}^{r_{1}}$ $...$ $%
e_{[i+m_{1}]}^{r_{1}}$ $e_{[i+m_{1}]}^{r_{2}}$ $...$ $%
e_{[i+m_{1}-m_{2}]}^{r_{2}}\}$ $\subset $ $D_{FP}(C_{N}^{\symbol{94}})$ is
computed by

\strut \strut 

$\ \ \ \ \ \ \ \ \ \ \mu _{\delta }\left(
\{e_{[i]}^{r_{1}}e_{[i+1]}^{r_{1}}...e_{[i+m_{1}]}^{r_{1}}e_{[i+m_{1}]}^{r_{2}}...e_{[i+m_{1}-m_{2}]}^{r_{2}}\}\right) =m_{1}+m_{2}.
$

\strut 

Notice that $\mu _{\delta }\left( \{\delta _{0}\}\right) $ $\leq $ $2N,$ for 
$\delta _{0}$ $\in $ $D_{FP}(C_{N}^{\symbol{94}}),$ different from that $\mu
_{C_{N}^{\symbol{94}}}(\{w\})$ $\in $ $\Bbb{N},$ for $w$ $\in $ $\Bbb{F}%
^{+}(C_{N}^{\symbol{94}}).$ Also, notice that the free semigroupoid $\Bbb{F}%
^{+}(C_{N}^{\symbol{94}})$ is an infinite set, but the diagram set $D(C_{N}^{%
\symbol{94}})$ is a finite set. Similar to the previous example, this
observation says that the energy measure $\mu _{C_{N}^{\symbol{94}}}$ of $%
C_{N}$ is not a bounded measure but the diagram measure $\mu _{\delta }$ is
a bounded measure. It is easy to check that $\mu _{C_{N}^{\symbol{94}}}$ is
locally bounded.
\end{example}

\strut \strut \strut \strut

\strut

\strut

\subsection{Graph Groupoid Measures and Reduced Diagram Measures}

\strut

\strut

\strut

Recall that the graph groupoid $\Bbb{G}$ of a finite directed graph $G$ is
the subset $\Bbb{F}_{r}^{+}(G^{\symbol{94}})$ of the free semigroupoid $\Bbb{%
F}^{+}(G^{\symbol{94}})$ of the shadowed graph $G^{\symbol{94}}$ of $G$
under the reducing relation (RR) with the inherited admissibility on $\Bbb{F}%
^{+}(G^{\symbol{94}}),$ and the reduced diagram set $D_{r}(G^{\symbol{94}})$
is the subset of the diagram set $D(G^{\symbol{94}})$ under (RR) with the
inherited admissibility on $D(G^{\symbol{94}})$. Notice that we can regard $%
\Bbb{G}$ and $D_{r}(G^{\symbol{94}})$ as subsets (sub-structures) of $\Bbb{F}%
^{+}(G^{\symbol{94}})$ and $D(G^{\symbol{94}}),$ respectively. So, we can
define the graph groupoid measure $\mu _{\Bbb{G}}$ and the reduced diagram
measure $\mu _{\delta ^{r}}$ by the restrictions of energy measure $\mu _{G^{%
\symbol{94}}}$ and the diagram measure $\mu _{\delta },$ respectively.

\strut

\begin{definition}
The graph groupoid measure $\mu _{\Bbb{G}}$ is defined by the restricted
measure $\mu _{G^{\symbol{94}}}\mid _{\Bbb{G}}$ of the energy measure $\mu
_{G^{\symbol{94}}}$ on the $\sigma $-algebra $P\left( \Bbb{G}\right) $. The
reduced diagram measure $\mu _{\delta ^{r}}$ is defined by the restricted
measure $\mu _{\delta }\mid _{D_{r}(G^{\symbol{94}})}$ of the diagram
measure $\mu _{\delta }$ on the $\sigma $-algebra $P\left( D_{r}(G^{\symbol{%
94}})\right) .$
\end{definition}

\strut \strut \strut

\begin{example}
Let $G_{N}$ be the one-vertex-$N$-loop-edge graph. i.e., the graph $G_{N}$
is the finite directed graph with $V(G_{N})$ $=$ $\{v\}$ and $E(G_{N})$ $=$ $%
\{l_{k}$ $=$ $v$ $l_{k}$ $v$ $:$ $k$ $=$ $1,$ ..., $N\}.$ We can have the
reduced finite path set

\strut 

$\ \ \ FP_{r}(G_{N}^{\symbol{94}})=\cup _{m=1}^{\infty }\left\{
l_{k_{1}}^{n_{1}}...l_{k_{m}}^{n_{m}}\left| 
\begin{array}{c}
k_{1}\neq k_{2},\,k_{2}\neq k_{3},...,\,k_{m-1}\neq k_{m} \\ 
k_{i}\in \{1,...,N\},\text{ \ \ }n_{i}\in \Bbb{Z} \\ 
i=1,\text{ }2,\text{ ..., }m,\text{ \ \ }m\in \Bbb{N}
\end{array}
\right. \right\} .$

\strut 

Thus the graph groupoid $\Bbb{G}_{N}$ of $G_{N}$ is $\{\emptyset \}$ $\cup $ 
$V(G_{N})$ $\cup $ $FP_{r}(G_{N}^{\symbol{94}}),$ as a set. We can easily
compute that the graph groupoid measure of $\{l_{k_{1}}^{n_{1}}$ ... $%
l_{k_{m}}^{n_{m}}\}$ $\subset $ $FP_{r}(G_{N}^{\symbol{94}})$ is

\strut 

$\ \ \ \ \ \ \ \ \ \ \ \ \ \ \ \ \ \ \ \ \ \ \mu _{\Bbb{G}_{N}}\left(
\{l_{k_{1}}^{n_{1}}...l_{k_{m}}^{n_{m}}\}\right) =\sum_{j=1}^{m}\left|
n_{j}\right| .$

\strut 

Let $D_{FP}^{r}(G_{N}^{\symbol{94}})$ be the subset $D_{r}(G_{N}^{\symbol{94}%
})$ $\setminus $ $(V(G_{N}^{\symbol{94}})$ $\cup $ $\{\emptyset \}).$ Then

\strut 

$\ \ \ D_{FP}^{r}(G_{N}^{\symbol{94}})=\cup _{m=1}^{N}\left\{
l_{k_{1}}^{r_{1}}...l_{k_{m}}^{r_{m}}\left| 
\begin{array}{c}
k_{1}\neq k_{2},\,k_{2}\neq k_{3},\,...,\,k_{m-1}\neq k_{m} \\ 
k_{i}\in \{1,...,N\},\text{ }r_{1},\text{...,}r_{m}\in \{\pm 1\} \\ 
\text{ }i=1,...,m,\text{ }m\in \Bbb{N}
\end{array}
\right. \right\} .$

\strut 

So, the reduced diagram measure of $\{l_{k_{1}}^{r_{1}}$ $...$ $%
l_{k_{m}}^{r_{m}}\}$ $\subset $ $D_{FP}^{r}(G_{N}^{\symbol{94}})$ is

\strut 

$\ \ \ \ \ \ \ \ \ \ \ \ \ \ \ \ \ \ \ \ \ \ \ \ \ \ \mu _{\delta
^{r}}\left( \{l_{k_{1}}^{r_{1}}...l_{k_{m}}^{r_{m}}\}\right) =m.$

\strut 

Remark that the reduced diagram $\delta ^{r}\left(
l_{k_{1}}^{n_{1}}...l_{k_{m}}^{n_{m}}\right) $ of $l_{k_{1}}^{n_{1}}$ ... $%
l_{k_{m}}^{n_{m}}$ $\in $ $\Bbb{G}$ is identical to the basic element $%
l_{k_{1}}^{r_{1}}$ ... $l_{k_{m}}^{r_{m}},$ contained in $D_{r}(G^{\symbol{94%
}}),$ whenever $k_{1}$ $\neq $ $k_{2},$ $k_{2}$ $\neq $ $k_{3},$ ..., $%
k_{m-1}$ $\neq $ $k_{m},$ and $r_{j}$ $=$ $1,$ if $n_{j}$ $>$ $0,$ and $r_{j}
$ $=$ $-1,$ if $n_{j}$ $<$ $0,$ for $j$ $=$ $1,$ ..., $m.$ Notice that the
graph groupoid $\Bbb{G}_{N}$ of $G_{N}$ is an infinite set, but the reduced
diagram set $D_{r}(G_{N}^{\symbol{94}})$ is a finite set. Hence the graph
groupoid measure $\mu _{\Bbb{G}_{N}}$ is not a bounded measure but the
diagram measure $\mu _{\delta ^{r}}$ is a bounded measure.
\end{example}

\strut

\begin{example}
Let $C_{N}$ be the one-flow circulant graph. The reduced finite path set of
the shadowed graph $C_{N}^{\symbol{94}}$ of $C_{N}$ is

\strut 

$\ \ \ \ \ \ FP_{r}(C_{N}^{\symbol{94}})=\cup _{m=0}^{\infty }\left\{ \left. 
\begin{array}{c}
e_{[i]}e_{[i+1]}\text{ ... }e_{[i+m]}, \\ 
e_{[i+m]}^{-1}e_{[i+m-1]}^{-1}...e_{[i]}^{-1}
\end{array}
\right. \left| 
\begin{array}{c}
i=1,...,N \\ 
\lbrack n]\in \Bbb{Z}_{N}
\end{array}
\right. \right\} .$

\strut 

Denote the graph groupoid of $C_{N}$ by $\Delta _{N}.$ Then the graph
groupoid measure of $\{e_{[i]}$ $e_{[i+1]}$ ... $e_{[i+m]}\}$ $\subset $ $%
FP_{r}(C_{N}^{\symbol{94}})$ is computed by

\strut 

$\ \ \ \ \ \ \ \ \ \ \ \ \ \ \ \ \ \ \ \ \ \mu _{\Delta _{N}}\left(
\{e_{[i]}e_{[i+1]}...e_{[i+m]}\}\right) $ $=$ $m.$

\strut 

The set $D_{FP}^{r}(C_{N}^{\symbol{94}})$ $\overset{def}{=}$ $D_{r}(C_{N}^{%
\symbol{94}})$ $\setminus $ $(V(C_{N}^{\symbol{94}})$ $\cup $ $\{\emptyset
\})$ is determined by

\strut 

$\ \ \ \ D_{FP}^{r}(C_{N}^{\symbol{94}})=\cup _{m=0}^{N-1}\left\{ \left. 
\begin{array}{c}
e_{[i]}e_{[i+1]}\text{ }...\text{ }e_{[i+m]}, \\ 
e_{[i+m]}^{-1}e_{[i+m-1]}^{-1}...e_{[i]}^{-1}
\end{array}
\right. \left| 
\begin{array}{c}
i=1,...,N \\ 
\lbrack n]\in \Bbb{Z}_{N}
\end{array}
\right. \right\} .$

\strut 

So, the diagram measure of $\{e_{[i]}$ $e_{[i+1]}$ ... $e_{[i+m]}\}$ $%
\subset $ $D_{FP}^{r}(C_{N}^{\symbol{94}})$ is computed by

\strut 

$\ \ \ \ \ \ \ \ \ \ \ \ \ \ \ \ \ \ \ \ \ \ \mu _{\delta ^{r}}\left(
\{e_{[i]}e_{[i+1]}...e_{[i+m]}\}\right) =m.$

\strut 

Notice that $\mu _{\delta ^{r}}\left( \{\delta _{0}\}\right) $ $\leq $ $N,$
for $\delta _{0}$ $\in $ $D_{FP}(C_{N}^{\symbol{94}}),$ different from that $%
\mu _{\Delta _{N}}(\{w\})$ $\in $ $\Bbb{N},$ for $w$ $\in $ $\Bbb{F}%
^{+}(C_{N}^{\symbol{94}}).$ The graph groupoid measure $\mu _{\Delta _{N}}$
of $C_{N}$ is not a bounded measure but the reduced diagram measure $\mu
_{\delta ^{r}}$ is a bounded measure.
\end{example}

\strut \strut \strut \strut

\strut

\strut

\subsection{Graph Measurings are Invariants on Graphs}

\strut

\strut

In this section, we will show that all graph measures we defined in the
previous sections are invariants on finite directed graphs. i.e., if we
denote $\mu _{G}$ as one of our graph measures, then $\mu _{G_{1}}$ and $\mu
_{G_{2}}$ are equivalent measures if and only if the graphs $G_{1}$ and $%
G_{2}$ are graph-isomorphic. Recall that two graphs $G_{1}$ and $G_{2}$ are
graph-isomorphic, if there is a graph-isomorphism $g$ $:$ $G_{1}$ $%
\rightarrow $ $G_{2}$ such that (i) $g$ is a bijection from $V(G_{1})$ onto $%
V(G^{\symbol{94}}),$ (ii) $g$ is also a bijection from $E(G_{1})$ onto $%
E(G_{2}),$ and (iii) if $e$ $=$ $v$ $e$ $v^{\prime }$ in $E(G_{1})$ with $v,$
$v^{\prime }$ $\in $ $V(G_{1}),$ then $g(e)$ $=$ $g(v)$ $g(e)$ $g(v^{\prime
})$ in $E(G_{2}).$

\strut

\begin{definition}
Let $G_{1}$ and $G_{2}$ be finite directed graphs and $G_{1}^{\symbol{94}},$ 
$G_{2}^{\symbol{94}}$, the corresponding shadowed graphs and let $\Bbb{F}%
^{+}(G_{k}^{\symbol{94}})$, $D(G_{k}^{\symbol{94}})$ and $\Bbb{G}_{k}$, $%
D_{r}(G_{k}^{\symbol{94}})$ be the free semigroupoid of $G_{k}^{\symbol{94}},
$ the diagram set of $G_{k}^{\symbol{94}}$ and the graph groupoid of $G_{k}$%
, the reduced diagram set of $G_{k},$ respectively, for $k$ $=$ $1,$ $2.$

\strut 

(1) We say that the free semigroupoids $\Bbb{F}^{+}(G_{k}^{\symbol{94}})$'s
are isomorphic if there exists a morphism $g^{\symbol{94}}$ $:$ $\Bbb{F}%
^{+}(G_{1}^{\symbol{94}})$ $\rightarrow $ $\Bbb{F}^{+}(G_{2}^{\symbol{94}})$
such that $g^{\symbol{94}}$ is bijective and $g^{\symbol{94}}$ preserves the
admissibility on $\Bbb{F}^{+}(G_{1}^{\symbol{94}})$ to that on $\Bbb{F}%
^{+}(G_{2}^{\symbol{94}}).$ i.e., $g^{\symbol{94}}(w_{1}$ $w_{2})$ $=$ $g^{%
\symbol{94}}(w_{1})$ $g^{\symbol{94}}(w_{2}),$ for all $w_{1},$ $w_{2}$ $\in 
$ $\Bbb{F}^{+}(G_{1}^{\symbol{94}}),$ with $g^{\symbol{94}}(\emptyset )$ $=$ 
$\emptyset .$

\strut 

(2) The diagram sets $D(G_{k}^{\symbol{94}})$'s are said to be isomorphic if
there exists a morphism $g_{\delta }$ $:$ $D(G_{1}^{\symbol{94}})$ $%
\rightarrow $ $D(G_{2}^{\symbol{94}})$ such that $g_{\delta }$ is bijective
and $g_{\delta }$ preserves the admissibility on $D(G_{1}^{\symbol{94}})$ to
that of $D(G_{2}^{\symbol{94}}).$ i.e., $g_{\delta }(\delta _{1}$ $\delta
_{2})$ $=$ $g_{\delta }(\delta _{1})$ $g_{\delta }(\delta _{2}),$ for all $%
\delta _{1}$ $\delta _{2}$ $\in $ $D(G_{1}^{\symbol{94}}).$

\strut 

(3) We say that the graph groupoids $\Bbb{G}_{k}$'s are isomorphic if there
exists a morphism $\Bbb{g}$ $:$ $\Bbb{G}_{1}$ $\rightarrow $ $\Bbb{G}_{2}$
such that $\Bbb{g}$ is bijective and it preserves the admissibility on $\Bbb{%
G}_{1}$ to that on $\Bbb{G}_{2}.$

\strut 

(4) The reduced diagram sets $D_{r}(G_{k}^{\symbol{94}})$'s are isomorphic
if there exists a morphism $g_{\delta ^{r}}$ $:$ $D_{r}(G_{1}^{\symbol{94}})$
$\rightarrow $ $D_{r}(G_{2}^{\symbol{94}})$ such that $g_{\delta ^{r}}$ is
bijective and it preserves the admissibility on $D_{r}(G_{1}^{\symbol{94}})$
to that on $D_{r}(G_{2}^{\symbol{94}}).$
\end{definition}

\strut

Let $X$ be an arbitrary set and assume that $X$ is partitioned by its subsets%
$\ X_{1}$ and $X_{2}.$ i.e., the set $X$ is the disjoint union of $X_{1}$
and $X_{2}.$ Let $f$ $:$ $X$ $\rightarrow $ $Y$ be a function, where $Y$ is
a set. Denote a function $f\ $by $f_{1}$ $\cup $ $f_{2}$ when $f$ satisfies
that

\strut

\begin{center}
$f(x)=\left( f_{1}\text{ }\cup \text{ }f_{2}\right) (x)\overset{def}{=}%
\left\{ 
\begin{array}{ll}
f_{1}(x) & \text{if }x\in X_{1} \\ 
f_{2}(x) & \text{if }x\in X_{2},
\end{array}
\right. $
\end{center}

\strut

for all $x$ $\in $ $X,$ and $f_{1}(X)$ $\cap $ $f_{2}(X)$ $=$ $\varnothing .$

\strut \strut \strut \strut \strut

\begin{theorem}
Let $G_{1}$ and $G_{2}$ be finite directed graphs. The shadowed graphs $%
G_{1}^{\symbol{94}}$ and $G_{2}^{\symbol{94}}$ are graph-isomorphic if and
only if the graph measures $\mu _{S_{1}}$ and $\mu _{S_{2}}$ are equivalent,
where $S_{k}$ $\in $ $\{G_{k}^{\symbol{94}},$ $\delta _{k},$ $\Bbb{G}_{k},$ $%
\delta _{k}^{r}\}$, for $k$ $=$ $1,$ $2.$
\end{theorem}

\strut

\begin{proof}
($\Rightarrow $) Suppose that the shadowed graph $G_{1}^{\symbol{94}}$ and $%
G_{2}^{\symbol{94}}$ are graph-isomorphic via the graph-isomorphism $g$ $:$ $%
G_{1}^{\symbol{94}}$ $\rightarrow $ $G_{2}^{\symbol{94}}.$ By this
graph-isomorphism $g$, we can construct a morphism $g^{\symbol{94}}$ $:$ $%
\Bbb{F}^{+}(G_{1}^{\symbol{94}})$ $\rightarrow $ $\Bbb{F}^{+}(G_{2}^{\symbol{%
94}})$ defined by

\strut

(1.5) $\ \ \ \ \ \ \ g^{\symbol{94}}(w)=\left\{ 
\begin{array}{ll}
g(w) & \text{if }w\in V(G_{1}^{\symbol{94}}) \\ 
g(w) & \text{if }w\in E(G_{1}^{\symbol{94}}) \\ 
g(e_{1})...g(e_{k}) & \text{if }w=e_{1}...e_{k}\in FP(G_{1}^{\symbol{94}})
\\ 
\emptyset & \text{if }w=\emptyset .
\end{array}
\right. $

\strut

Then, by the bijectivity of $g$ and by the admissibility-preserving property
on of $g,$ the map $g^{\symbol{94}}$ is admissibility-preserving bijection,
and hence $\Bbb{F}^{+}(G_{1}^{\symbol{94}})$ and $\Bbb{F}^{+}(G_{2}^{\symbol{%
94}})$ are isomorphic. Therefore, there exists a map $\Phi ^{\symbol{94}}$ $%
: $ $P\left( \Bbb{F}^{+}(G_{1}^{\symbol{94}})\right) $ $\rightarrow $ $%
P\left( \Bbb{F}^{+}(G_{2}^{\symbol{94}})\right) $ from the $\sigma $-algebra 
$P\left( \Bbb{F}^{+}(G_{1}^{\symbol{94}})\right) $ onto the $\sigma $%
-algebra $P\left( \Bbb{F}^{+}(G_{2}^{\symbol{94}})\right) $ defined by

\strut

(1.6) \ \ \ \ \ \ \ $\ \Phi ^{\symbol{94}}\left( S\right) \overset{def}{=}%
\{g^{\symbol{94}}(w):w\in S\},$ for all $S$ $\in $ $P\left( \Bbb{F}%
^{+}(G_{1}^{\symbol{94}})\right) ,$

\strut

such that

\strut

$\ \ \ \ \ \ \ \ \ \ \ 
\begin{array}{ll}
\mu _{G_{2}^{\symbol{94}}}\left( \Phi ^{\symbol{94}}(S)\right) & =\mu
_{G_{2}^{\symbol{94}}}\left( \{g^{\symbol{94}}(w):w\in S\}\right) \\ 
&  \\ 
& =\underset{w\in S}{\sum }\mu _{G_{2}^{\symbol{94}}}\left( \{g^{\symbol{94}%
}(w)\}\right) =\underset{w\in S}{\sum }\mu _{G_{1}^{\symbol{94}}}(\{w\}),
\end{array}
$

\strut

by (1.5). Clearly, by the bijectivity of $g^{\symbol{94}},$ the map $\Phi ^{%
\symbol{94}}$ is also bijective. Therefore, the energy measures $\mu
_{G_{1}^{\symbol{94}}}$ and $\mu _{G_{2}^{\symbol{94}}}$ are equivalent.

\strut

Since the graph groupoid $\Bbb{G}_{k}$ are the structures $\Bbb{F}%
_{r}^{+}(G_{k}^{\symbol{94}})$ (with reducing relation (RR)), the graph
groupoid measures $\mu _{\Bbb{G}_{k}}$'s are equivalent similarly, for $k$ $%
= $ $1,$ $2.$

\strut

Also, if $G_{1}$ and $G_{2}$ are graph-isomorphic, then we can determine the
morphism $g_{\delta }$ $:$ $D(G_{1}^{\symbol{94}})$ $\rightarrow $ $D(G_{2}^{%
\symbol{94}})$ defined by

\strut

(1.7) $\ \ \ \ \ \ \ \ \ \ \ \ \ \ \ g_{\delta }$ $(\delta _{w}^{1})$ $=$ $%
\delta _{g^{\symbol{94}}(w)}^{2},$ for all $\delta _{w}^{1}$ $\in $ $%
D(G_{1}^{\symbol{94}}).$

\strut

As we observed in Section 1.1, diagram sets are quotient structure of free
semigroupoids. So, any diagrams in the diagram set are regarded as
equivalence classes of elements in a free semigroupoid having the same
diagrams. So, if $\delta _{0}^{1}$ is an element of $D(G_{1}^{\symbol{94}}),$
then there always exists $w$ $\in $ $\Bbb{F}^{+}(G_{1}^{\symbol{94}})$ such
that $\delta _{0}^{1}$ $=$ $\delta _{w}^{1}.$ By the morphism $g^{\symbol{94}%
}$ defined in (1.5), we can determine $g^{\symbol{94}}(\delta _{w}^{1})$ $%
\subseteq $ $D(G_{2}^{\symbol{94}}).$ So, the map $g_{\delta }$ is a
well-defined bijection between $D(G_{1}^{\symbol{94}})$ and $D(G_{2}^{%
\symbol{94}}),$ and moreover it is a admissibility-preserving morphisms.
i.e.,

\strut

$\ \ \ \ \ \ \ \ \ \ \ \ 
\begin{array}{ll}
g_{\delta }(\delta _{w_{1}}^{1}\delta _{w_{2}}^{1}) & =g_{\delta }(\delta
_{w_{1}w_{2}}^{1})=\delta _{g^{\symbol{94}}(w_{1}w_{2})}^{2}=\delta _{g^{%
\symbol{94}}(w_{1})g^{\symbol{94}}(w_{2})}^{2} \\ 
&  \\ 
& =\delta _{g^{\symbol{94}}(w_{1})}^{2}\delta _{g^{\symbol{94}%
}(w_{2})}^{2}=g_{\delta }(\delta _{w_{1}}^{1})g_{\delta }(\delta
_{w_{2}}^{1}).
\end{array}
$

$\strut $

Therefore, the diagram sets $D(G_{1}^{\symbol{94}})$ and $D(G_{2}^{\symbol{94%
}})$ are isomorphic by the morphism $g_{\delta }.$ So, we can define a map $%
\Phi _{\delta }$ from the $\sigma $-algebra $P\left( D(G_{1}^{\symbol{94}%
})\right) $ onto the $\sigma $-algebra $P\left( D(G_{2}^{\symbol{94}%
})\right) $ by

\strut

(1.8) \ \ \ \ \ \ \ \ $\ \Phi _{\delta }(S)=\{g_{\delta }(\delta
_{0}^{1}):\delta _{0}^{1}\in S\},$ for all $S$ $\in $ $D(G_{1}^{\symbol{94}%
}).$

\strut

Then, by $g_{\delta },$ the map $\Phi _{\delta }$ is bijective. Also we can
have that

\strut

$\ 
\begin{array}{ll}
\mu _{\delta ^{2}}\left( \Phi _{\delta }(S)\right) & =\mu _{\delta
^{2}}\left( \{g_{\delta }(\delta _{0}^{1}):\delta _{0}^{1}\in S\}\right) =%
\underset{\delta _{0}^{1}\in S}{\sum }\mu _{\delta ^{2}}\left( \{g_{\delta
}(\delta _{0}^{1})\}\right) \\ 
&  \\ 
& =\underset{\delta _{w}^{1}\in S}{\sum }\mu _{\delta ^{2}}\left( \{\delta
_{g_{\delta }(w)}^{2}\}\right) =\underset{\delta _{w}^{1}\in S}{\sum }\mu
_{\delta ^{1}}\left( \{\delta _{w}^{1}\}\right) =\mu _{\delta ^{1}}\left(
S\right) .
\end{array}
$

\strut

Thus the measures $\mu _{\delta ^{1}}$ and $\mu _{\delta ^{2}}$ are
equivalent. Under the reducing relation (RR), similarly, we can conclude
that the reduced diagram measures $\mu _{\delta ^{r:k}}$ on $P\left(
D_{r}(G_{k}^{\symbol{94}})\right) ,$ for $k$ $=$ $1,$ $2,$ are equivalent.

\strut

($\Leftarrow $) Assume now that the graph measures $\mu _{G_{1}^{\symbol{94}%
}}$ and $\mu _{G_{2}^{\symbol{94}}}$ are equivalent. i.e., there exists a
bijection $\Phi ^{\symbol{94}}$ $:$ $P\left( \Bbb{F}^{+}(G_{1}^{\symbol{94}%
})\right) $ $\rightarrow $ $P\left( \Bbb{F}^{+}(G_{2}^{\symbol{94}})\right) $
such that

\strut

(1.9) \ $\ \ \ \ \ \ \ \mu _{G_{2}^{\symbol{94}}}\left( \Phi ^{\symbol{94}%
}(S)\right) $ $=$ $\mu _{G_{1}^{\symbol{94}}}\left( S\right) ,$ for all $S$ $%
\in $ $P\left( \Bbb{F}^{+}(G_{1}^{\symbol{94}})\right) .$

\strut

Then since

\strut \strut

$P\left( \Bbb{F}^{+}(G_{k}^{\symbol{94}})\right) =\{S_{V}\cup
S_{FP}:S_{V}\in P\left( V(G_{k}^{\symbol{94}})\right) ,$ $S_{FP}\in P\left(
FP(G_{k}^{\symbol{94}})\right) \}$,

\strut

and

\strut

$\ \ \ \ \ \ \ \ \ \ \ \ \ \ \ \ \ \ \ \ \ \ \ P\left( V(G_{k}^{\symbol{94}%
})\right) \cap P\left( FP(G_{k}^{\symbol{94}})\right) =\varnothing ,$

\strut

for $k$ $=$ $1,$ $2,$ the bijection $\Phi ^{\symbol{94}}$ is re-written as $%
\Phi ^{\symbol{94}}$ $=$ $\Phi _{V}^{\symbol{94}}$ $\cup $ $\Phi _{FP}^{%
\symbol{94}},$ where $\Phi _{V}^{\symbol{94}}$ $=$ $\Phi ^{\symbol{94}}$ $%
\mid _{P\left( V(G_{1}^{\symbol{94}})\right) }$ and $\Phi _{FP}^{\symbol{94}%
} $ $=$ $\Phi ^{\symbol{94}}$ $\mid _{P\left( FP(G_{1}^{\symbol{94}})\right)
}, $ defined by

\strut

$\ \Phi ^{\symbol{94}}\left( S_{V}\cup S_{FP}\right) =\left( \Phi _{V}^{%
\symbol{94}}\cup \Phi _{FP}^{\symbol{94}}\right) \left( S_{V}\cup
S_{FP}\right) \overset{def}{=}\Phi _{V}^{\symbol{94}}(S_{V})+\Phi _{FP}^{%
\symbol{94}}(S_{FP}),$

\strut \strut

for all $S_{V}\cup S_{FP}$ $\in $ $P\left( \Bbb{F}^{+}(G_{1}^{\symbol{94}%
})\right) .$ With respect to $\Phi _{V}^{\symbol{94}},$ we can construct the
bijection $g_{V}^{\symbol{94}}$ $:$ $V(G_{1}^{\symbol{94}})$ $\rightarrow $ $%
V(G_{2}^{\symbol{94}})$ defined by

\strut

$\ \ \ \ \ \ \ \ \ \ \ \ \ \ \ g_{V}^{\symbol{94}}$ $(v_{1})$ $\overset{def}{%
=}$ $\Phi _{V}^{\symbol{94}}(\{v_{1}\}),$ $\ $\ for all $\ \ v_{1}$ $\in $ $%
V(G_{1}^{\symbol{94}}),$

\strut

and similarly, with respect to $\Phi _{FP}^{\symbol{94}},$ we can construct
the bijection $g_{FP}^{\symbol{94}}$ $:$ $FP(G_{1}^{\symbol{94}})$ $%
\rightarrow $ $FP(G_{2}^{\symbol{94}})$ defined by

\strut

$\ \ \ \ \ \ \ \ \ \ \ \ \ g_{FP}^{\symbol{94}}(w_{1})$ $\overset{def}{=}$ $%
\Phi _{FP}^{\symbol{94}}(\{w_{1}\}),$ for all $w_{1}\in FP(G_{1}^{\symbol{94}%
}).$

\strut

Then we have a bijection $g^{\symbol{94}}$ $=$ $g_{V}^{\symbol{94}}$ $\cup $ 
$g_{FP}^{\symbol{94}}$ from $\Bbb{F}^{+}(G_{1}^{\symbol{94}})$ onto $\Bbb{F}%
^{+}(G_{2}^{\symbol{94}}).$ Let's suppose that $g^{\symbol{94}}$ does not
preserve the admissibility. Then clearly $\Phi ^{\symbol{94}}$ does not
satisfy (1.9), when we take $S$ $=$ $\Bbb{F}^{+}(G_{1}^{\symbol{94}}).$ This
contradict our assumption that $\mu _{G_{k}^{\symbol{94}}}$'s are
equivalent, for $k$ $=$ $1,$ $2.$ So, the map $g^{\symbol{94}}$ is the
admissibility-preserving bijective morphism. From this morphism $g^{\symbol{%
94}},$ we can construct a graph-isomorphism $g$ $:$ $G_{1}^{\symbol{94}}$ $%
\rightarrow $ $G_{2}^{\symbol{94}}$ such that $g$ $=$ $g^{\symbol{94}}$ $%
\mid _{V(G^{\symbol{94}})\cup E(G^{\symbol{94}})}.$

\strut

Similar to the previous proof, we can conclude that if the graph groupoid
measures $\mu _{\Bbb{G}_{1}}$and $\mu _{\Bbb{G}_{2}}$ are equivalent, then
the shadowed graphs $G_{1}^{\symbol{94}}$ and $G_{2}^{\symbol{94}}$ are
graph-isomorphic.

\strut

Suppose now that the diagram measures $\mu _{\delta ^{1}}$ and $\mu _{\delta
^{2}}$ are equivalent. i.e., there exists a bijection $\Phi _{\delta }$ from
the $\sigma $-algebra $P\left( D(G_{1}^{\symbol{94}})\right) $ onto the $%
\sigma $-algebra $P\left( D(G_{2}^{\symbol{94}})\right) $ satisfying that

\strut

(1.9)$^{\prime }$ $\ \ \ \ \ \ \mu _{\delta ^{2}}$ $(\Phi _{\delta }(S))$ $=$
$\mu _{\delta ^{1}}$ $(S),$ for all $S$ $\in $ $P\left( D(G_{1}^{\symbol{94}%
})\right) .$

\strut

So, we can construct a map $g_{\delta }$ $:$ $D(G_{1}^{\symbol{94}})$ $%
\rightarrow $ $D(G_{2}^{\symbol{94}})$ defined by

\strut

$\ \ \ \ \ \ \ \ \ g_{\delta }(\delta _{1})$ $=$ $\delta _{2}\in D(G_{2}^{%
\symbol{94}}),$ whenever $\Phi _{\delta }$ $(\{\delta _{1}\})=\{\delta
_{2}\},$

\strut

for all $\delta _{1}$ $\in $ $D(G_{1}^{\symbol{94}}).$ Since $\Phi _{\delta
} $ is bijective from $\{\{\delta _{1}\}$ $:$ $\delta _{1}$ $\in $ $D(G_{1}^{%
\symbol{94}})\}$ onto $\{\{\delta _{2}\}$ $:$ $\delta _{2}$ $\in $ $D(G_{2}^{%
\symbol{94}})\},$ the map $g_{\delta }$ is bijective, too. Assume now that $%
g_{\delta }$ does not preserves the admissibility. Then this contradict (1.9)%
$^{\prime }.$ Also since $D(G_{k}^{\symbol{94}})$ is the disjoint union of $%
V(G_{k}^{\symbol{94}})$ and $D_{FP}(G_{k}^{\symbol{94}}),$ where $%
D_{FP}(G_{k}^{\symbol{94}})$ $=$ $D(G_{k}^{\symbol{94}})$ $\cap $ $FP(G_{k}^{%
\symbol{94}}),$ there are $g_{\delta :V}$ $=$ $g_{\delta }$ $\mid _{V(G_{1}^{%
\symbol{94}})}$ and $g_{\delta :FP}$ $=$ $g_{\delta }$ $\mid _{D_{FP}(G_{1}^{%
\symbol{94}})}$ such that $g_{\delta }$ $=$ $g_{\delta :V}$ $\cup $ $%
g_{\delta :FP}.$ Define $g^{\symbol{94}}$ $:$ $G_{1}^{\symbol{94}}$ $%
\rightarrow $ $G_{2}^{\symbol{94}}$ by $g^{\symbol{94}}$ $=$ $g_{\delta :V}$ 
$\cup $ $g_{\delta :E^{\symbol{94}}},$ where $g_{\delta :E^{\symbol{94}}}$ $%
= $ $g_{\delta :FP}$ $\mid _{E(G_{1}^{\symbol{94}})}.$ Then it is a
graph-isomorphism from $G_{1}^{\symbol{94}}$ onto $G_{2}^{\symbol{94}}.$

\strut

Similarly, if the reduced diagram measures $\mu _{\delta _{1}^{r}}$ and $\mu
_{\delta _{2}^{r}}$ of $G_{1}$ and $G_{2}$ are equivalent, then the shadowed
graphs $G_{1}^{\symbol{94}}$ and $G_{2}^{\symbol{94}}$ are graph-isomorphic.
\end{proof}

\strut \strut

The above theorem shows that our graph measurings are invariants on shadowed
graphs of finite directed graphs. Also, the proof of the previous theorem
provides the following corollary.

\strut

\begin{corollary}
Let $G_{1}$ and $G_{2}$ be finite directed graphs. Then the shadowed graphs $%
G_{1}^{\symbol{94}}$ and $G_{2}^{\symbol{94}}$ are graph-isomorphic if and
only if the free semigroupoids $\Bbb{F}^{+}(G_{1}^{\symbol{94}})$ and $\Bbb{F%
}^{+}(G_{2}^{\symbol{94}})$ are isomorphic if and only if the graph
groupoids $\Bbb{G}_{1}$ and $\Bbb{G}_{2}$ are isomorphic if and only if the
diagram sets $D(G_{1}^{\symbol{94}})$ and $D(G_{2}^{\symbol{94}})$ are
isomorphic if and only if the reduced diagram sets $D_{r}(G_{1}^{\symbol{94}%
})$ and $D_{r}(G_{2}^{\symbol{94}})$ are isomorphic. $\square $
\end{corollary}

\strut

The corollary shows that the constructions of algebraic structures
introduced in this paper are all invariants on \emph{shadowed graphs} of
finite directed graphs. It is easy to check that if two finite directed
graphs $G_{1}$ and $G_{2}$ are graph-isomorphic, then the shadows $%
G_{1}^{-1} $ and $G_{2}^{-1}$ are graph-isomorphic, and hence the shadowed
graphs $G_{1}^{\symbol{94}}$ and $G_{2}^{\symbol{94}}$ are graph-isomorphic.

\strut

Suppose the graphs $G_{1}$ and $G_{2}$ are graph-isomorphic by the
graph-isomorphism $g_{0}^{(+1)}$ $:$ $G_{1}$ $\rightarrow $ $G_{2}.$ Then we
have the graph-isomorphism $g_{0}^{(-1)}$ from the shadow $G_{1}^{-1}$ of $%
G_{1}$ onto the shadow $G_{2}^{-1}$ of $G_{2}.$ Since the admissibility on $%
G_{k}^{-1}$ are oppositely preserved by that of $G_{k},$ for $k$ $=$ $1,$ $%
2, $ we can naturally define the isomorphism $g_{0}^{(-1)}$ by

\strut

$\ \ \ \ \ \ \ \ \ \ \ g_{0}^{(-1)}(v)$ $\overset{def}{=}$ $g_{0}^{(+1)}(v),$
\ for all \ $v$ $\in $ $V(G_{1}^{-1})$ $=$ $V(G_{1})$

and

$\ \ \ \ \ \ \ \ \ \ \ g_{0}^{(-1)}(e^{-1})\overset{def}{=}\left(
g_{0}^{(+1)}(e)\right) ^{-1},$ for all $e$ $\in $ $E(G_{1}^{-1}).$

\strut

Notice that we can regard $g_{0}^{(\pm 1)}$ by $g_{0}^{(\pm 1)}$ $=$ $%
g_{0:V(G_{1}^{\pm 1})}^{(\pm 1)}$ $\cup $ $g_{0:E(G_{1}^{\pm 1})}^{(\pm 1)},$
where $g_{0:V(G_{1})}^{(\pm 1)}$ $=$ $g_{0}^{(\pm 1)}$ $\mid _{V(G_{1}^{\pm
1})}$ and $g_{0:E(G_{1}^{\pm 1})}^{(\pm 1)}$ $=$ $g_{0}^{(\pm 1)}$ $\mid
_{E(G_{1}^{\pm 1})},$ since $V(G_{k})$ and $E(G_{k})$ are disjoint, for $k$ $%
=$ $1,$ $2.$ Then there exists a graph-isomorphism $g$ $:$ $G_{1}^{\symbol{94%
}}$ $\rightarrow $ $G_{2}^{\symbol{94}}$ defined by $g$ $=$ $%
g_{0:V(G_{1})}^{(+1)}$ $\cup $ $g_{0:E(G_{1})}^{(+1)}$ $\cup $ $%
g_{0:E(G_{1}^{-1})}^{(-1)},$ since $E(G_{k})$ and $E(G_{k}^{-1})$ are
disjoint, for $k$ $=$ $1,$ $2.$ So, indeed, the shadowed graphs $G_{1}^{%
\symbol{94}}$ and $G_{2}^{\symbol{94}}$ are graph-isomorphic, whenever $%
G_{1} $ and $G_{2}$ are graph-isomorphic. However, the converse does not
hold true, in general.

\strut

\begin{example}
Let $G_{1}$ be a one-flow circulant graph $C_{3}$ with $V(G_{1})$ $=$ $%
\{v_{1},$ $v_{2},$ $v_{3}\}$ and $E(G_{1})$ $=$ $\{e_{1}$ $=$ $v_{1}$ $e_{1}$
$v_{2},$ $e_{2}$ $=$ $v_{2}$ $e_{2}$ $v_{3},$ $e_{3}$ $=$ $v_{3}$ $e_{3}$ $%
v_{1}\}.$ And let $G_{2}$ be a finite directed graph with $V(G_{2})$ $=$ $%
\{v_{1}^{\prime },$ $v_{2}^{\prime },$ $v_{3}^{\prime }\}$ and $E(G_{2})$ $=$
$\{e_{1}^{\prime }$ $=$ $v_{1}^{\prime }$ $e_{1}^{\prime }$ $v_{2}^{\prime },
$ $e_{2}^{\prime }$ $=$ $v_{2}^{\prime }$ $e_{1}^{\prime }$ $v_{1}^{\prime },
$ $e_{3}^{\prime }$ $=$ $v_{3}^{\prime }$ $e_{3}^{\prime }$ $v_{1}^{\prime
}\}.$ Then these two graphs are not graph-isomorphic. However, the shadowed
graphs $G_{1}^{\symbol{94}}$ and $G_{2}^{\symbol{94}}$ are graph-isomorphic,
since there exists a graph-isomorphism $g$ $:$ $G_{1}^{\symbol{94}}$ $%
\rightarrow $ $G_{2}^{\symbol{94}}$ such that

\strut 

\ $\ \ \ \ \ \ \ \ \ \ \ \ \ \ \ \ \ \ \ g(v_{k})$ $=$ $v_{k}^{\prime },$ \
for all \ $k$ $=$ $1,$ $2,$ $3$

and

$\ \ \ \ \ \ \ g(e_{1}^{\pm 1})$ $=$ $(e_{1}^{\prime })^{\pm 1},$ \ $%
g(e_{2}^{\pm 1})$ $=$ $(e_{2}^{\prime })^{\mp 1}$ and $g(e_{3}^{\pm 1})$ $=$ 
$(e_{3}^{\prime })^{\pm 1}.$

\strut 

So, even though $G_{1}$ and $G_{2}$ are not graph-isomorphic, the shadowed
graphs $G_{1}^{\symbol{94}}$ and $G_{2}^{\symbol{94}}$ are graph-isomorphic.
This shows that we cannot guarantee that if $G_{1}^{\symbol{94}}$ and $%
G_{2}^{\symbol{94}}$ are graph-isomorphic, then $G_{1}$ and $G_{2}$ are
graph-isomorphic.
\end{example}

\strut

Therefore, we cannot conclude that our graph measurings are invariants on
finite directed graphs. However, if we consider shadowed graphs as
two-colored graphs, then we can get that our graph measurings are invariants
on finite directed graphs. i.e., we can show that if two shadowed graphs $%
G_{1}^{\symbol{94}}$ and $G_{2}^{\symbol{94}}$ are
two-colored-graph-isomorphic (under certain coloring), then the graphs $G_{1}
$ and $G_{2}$ are graph-isomorphic.

\strut

\begin{definition}
Let $G$ be a (finite) directed graph. We say that the graph $G$ is colored
if each edge $e$ $\in $ $E(G)$ has its color $c_{e}.$ We assume that the
vertices of colored graphs are not colored. Let $G_{1}$ and $G_{2}$ be
finite directed colored graphs. Then the graphs $G_{1}$ and $G_{2}$ are
colored-graph-isomorphic if (i) $G_{1}$ and $G_{2}$ are graph-isomorphic via
the graph-isomorphism $g$ $:$ $G_{1}$ $\rightarrow $ $G_{2}$ and (ii) $%
g(c_{e})$ $=$ $c_{g(e)},$ for all $e$ $\in $ $E(G_{1}).$ The morphism $g$
satisfying (i) and (ii) is called the colored-graph-isomorphism.
\end{definition}

\strut

By definition, the colored-graph-isomorphisms are graph-isomorphisms
preserving the colorings. So, if two finite directed graphs $G_{1}$ and $%
G_{2}$ are colored-graph-isomorphic, then they are graph-isomorphic. Of
course, the converse does not hold. From now, we will assume that all finite
directed graphs are one-colored graphs having the same colorings. we will
regard shadowed graphs of finite directed graphs as two-colored graphs.

\strut

\begin{definition}
Let $G$ be a finite directed graph with its coloring $c_{1}$ and let $G^{-1}$
be the shadow of $G.$ The shadowed graph $G^{\symbol{94}}$ is a finite
directed graph having

\strut 

$\ \ \ \ \ \ V(G^{\symbol{94}})$ $=$ $V(G)$ $=$ $V(G^{-1})$ and $E(G^{%
\symbol{94}})$ $=$ $E(G)$ $\cup $ $E(G^{-1})$

\strut 

with the two-colors $c_{1}$ and $c_{2}$ such that

\strut 

$\ \ \ \ \ \ \ c_{1}=c_{e}$ \ and \ $c_{2}=c_{e^{-1}},$ \ for all \ $e$ $\in 
$ $E(G)$ $\subset $ $E(G^{\symbol{94}}),$

\strut 

where $c_{f}$ means the color of an edge $f$ $\in $ $E(G^{\symbol{94}}).$%
\strut 
\end{definition}

\strut

\begin{example}
Let $G_{1}$ be the one-flow circulant graph $C_{3}$ and $G_{2},$ a finite
directed graph given in the previous example. We do know that $G_{1}$ and $%
G_{2}$ are not graph-isomorphic. In the previous example, we could construct
a graph-isomorphism $g$ $:$ $G_{1}^{\symbol{94}}$ $\rightarrow $ $G_{2}^{%
\symbol{94}},$ where $G_{k}^{\symbol{94}}$ are shadowed graphs of $G_{k},$
for $k$ $=$ $1,$ $2,$ without considering the shadowed graphs are colored.
Now, we consider $G_{1}^{\symbol{94}}$ and $G_{2}^{\symbol{94}}$ are
two-colored graphs having their colors $\{c_{1}^{(1)},$ $c_{2}^{(1)}\}$ and $%
\{c_{1}^{(2)},$ $c_{2}^{(2)}\},$ respectively. In particular, let $%
c_{1}^{(k)}$ be the color of all edges in $E(G_{k})$ and let $c_{2}^{(k)}$
be the color of all edges in $E(G_{k}^{-1}),$ for $k$ $=$ $1,$ $2.$ If $g$
were a colored-graph-isomorphism from $G_{1}^{\symbol{94}}$ to $G_{2}^{%
\symbol{94}},$ then $g$ should satisfy $g(c_{i}^{(1)})$ $=$ $c_{i}^{(2)},$
for $i$ $=$ $1,$ $2$. We can easily check that the graph-isomorphism $g$ in
the previous example does not preserve the coloring. Especially,

\strut 

$\ \ \ \ \ \ \ \ \ \ \ \ \ g(c_{e_{2}})$ $=$ $g(c_{1}^{(1)})$ $=$ $%
c_{2}^{(2)}$ $\neq $ $c_{1}^{(2)}$ $=$ $c_{g(e_{2})}$ $=$ $c_{e_{2}^{\prime
}},$

\strut 

for the edge $e_{2}$ in $E(G_{1}).$ So, by regarding shadowed graphs as
two-colored graphs, the map $g$ is not a colored-graph-isomorphism. Indeed,
we can see that there is no morphisms preserving the coloring and the
admissibility at the same time. In other words, the shadowed graphs $G_{1}^{%
\symbol{94}}$ and $G_{2}^{\symbol{94}}$ are not colored-graph-isomorphic and
hence they are not isomorphic as shadowed graphs in the sense of above new
definition.
\end{example}

\strut \strut \strut

In the rest of this section, we will assume shadowed graphs are two-colored
graphs such that all edges of the given graphs have the same color $c_{1}$
and the shadows have the same color $c_{2}$, where $c_{1}$ $\neq $ $c_{2}.$

\strut

\begin{proposition}
Let $G_{1}$ and $G_{2}$ be finite directed graphs and let $G_{1}^{\symbol{94}%
}$ and $G_{2}^{\symbol{94}}$ be the corresponding (two-colored) shadowed
graphs having the same two-coloring. If the shadowed graphs $G_{1}^{\symbol{%
94}}$ and $G_{2}^{\symbol{94}}$ are colored-graph-isomorphic, then the
graphs $G_{1}$ and $G_{2}$ are graph-isomorphic.
\end{proposition}

\strut

\begin{proof}
Suppose $g^{\symbol{94}}$ $:$ $G_{1}^{\symbol{94}}$ $\rightarrow $ $G_{2}^{%
\symbol{94}}$ is a colored-graph-isomorphism. Since $V(G_{k}^{\symbol{94}})$
and $E(G_{k}^{\symbol{94}})$ are disjoint, for $k$ $=$ $1,$ $2,$ we can
rewrite $g^{\symbol{94}}$ by $g^{\symbol{94}}$ $=$ $g_{V}^{\symbol{94}}$ $%
\cup $ $g_{E}^{\symbol{94}},$ where $g_{V}^{\symbol{94}}$ $=$ $g^{\symbol{94}%
}$ $\mid _{V(G_{1}^{\symbol{94}})}$ and $g_{FP}^{\symbol{94}}$ $=$ $g^{%
\symbol{94}}$ $\mid _{E(G_{1}^{\symbol{94}})}.$ Notice that this morphism $%
g^{\symbol{94}}$ preserves the coloring of $G_{1}^{\symbol{94}}$ to that of $%
G_{2}^{\symbol{94}}.$ Since $E(G_{k}^{\symbol{94}})$ are partitioned $%
E(G_{k})$ and $E(G_{k}^{-1})$ and since $E(G_{k})$ and $E(G_{k}^{-1})$ have
different colors, for $k$ $=$ $1,$ $2,$ we can rewrite $g_{E}^{\symbol{94}}$
by $g_{E}^{\symbol{94}}$ $=$ $g_{+1}^{\symbol{94}}$ $\cup $ $g_{-1}^{\symbol{%
94}},$ where $g_{+}^{\symbol{94}}$ $=$ $g_{E}^{\symbol{94}}$ $\mid
_{E(G_{1})}$ and $g_{-}^{\symbol{94}}$ $=$ $g_{E}^{\symbol{94}}$ $\mid
_{E(G_{1}^{-1})}.$ Moreover, if $G_{k}$ has two colors $c_{1}^{(k)}$ and $%
c_{2}^{(k)},$ where $c_{1}^{(k)}$ $=$ $c_{e_{k}}$ and $c_{2}^{(k)}$ $=$ $%
c_{e_{k}^{-1}},$ for all $e_{k}$ $\in $ $E(G_{k})$, for $k$ $=$ $1,$ $2,$
then we can have the map $g_{+1}^{\symbol{94}}$ and $g_{-1}^{\symbol{94}}$
satisfy that $g_{+1}^{\symbol{94}}(c_{1}^{(1)})$ $=$ $c_{1}^{(2)}$ and $%
g_{-1}^{\symbol{94}}(c_{2}^{(1)})$ $=$ $c_{2}^{(2)}.$ Now, we can define the
morphism $g$ $:$ $G_{1}$ $\rightarrow $ $G_{2}$ by $g$ $=$ $g_{V}^{\symbol{94%
}}$ $\cup $ $g_{+1}^{\symbol{94}}.$ Then this map $g$ is a graph-isomorphism.
\end{proof}

\strut

\begin{corollary}
Let $G_{1}$ and $G_{2}$ be finite directed graphs. Then the graphs $G_{1}$
and $G_{2}$ are graph-isomorphic if and only if the shadowed graphs $G_{1}^{%
\symbol{94}}$ and $G_{2}^{\symbol{94}}$ are colored-graph-isomorphic. $%
\square $
\end{corollary}

\strut \strut \strut

By the previous theorem and by the previous proposition, if we regard
shadowed graphs as two-colored graphs, then our graph measurings are
invariants on finite directed graphs.

\strut

\begin{theorem}
Let $G_{1}$ and $G_{2}$ be finite directed graphs. Then the graphs $G_{1}$
and $G_{2}$ are graph-isomorphic if and only if the graph measures $\mu
_{S_{1}}$ and $\mu _{S_{2}}$ are equivalent, where $S_{k}$ $\in $ $\{G_{k}^{%
\symbol{94}},$ $\delta _{k},$ $\Bbb{G}_{k},$ $\delta _{k}^{r}\},$ for $k$ $=$
$1,$ $2.$
\end{theorem}

\strut

\begin{proof}
By the previous proposition, two finite directed graphs $G_{1}$ and $G_{2}$
are graph-isomorphic if and only if two (two-colored) shadowed graphs $%
G_{1}^{\symbol{94}}$ and $G_{2}^{\symbol{94}}$ are colored-graph-isomorphic.
By the previous theorem, the shadowed graphs $G_{1}^{\symbol{94}}$ and $%
G_{2}^{\symbol{94}}$ are (colored) graph-isomorphic if and only if the graph
measures $\mu _{S_{1}}$ and $\mu _{S_{2}}$ are equivalent.
\end{proof}

\strut \strut \strut \strut

\begin{remark}
The two-colored graph setting of shadowed graphs provides us that a shadowed
graph $G^{\symbol{94}}$ of $G$ is uniquely determined up to
graph-isomorphisms.
\end{remark}

\strut \strut

\strut

\strut

\section{Reduced Diagram Measure Theory on Graphs}

\strut

\strut

In this chapter, we will consider Reduced Diagram Measure Theory on a finite
directed graph $G$, with respect to the reduced diagram measure $\mu
_{\delta ^{r}}$ of $G.$ Other graph measures will have same or similar
Measure Theory. For convenience, we will denote the reduced diagram measure $%
\mu _{\delta ^{r}}$ by $\mu _{r}$ and we will denote the measure space $%
(D_{r}(G^{\symbol{94}}),$ $P\left( D_{r}(G^{\symbol{94}})\right) ,$ $\mu
_{r})$ by $(G^{\symbol{94}},$ $\mu _{r}).$ If other graph measures have
different results, then we mention about them in remarks. The reason why we
concentrate on observing Measure Theory with respect to reduced diagram
measures is that reduced diagram measures are bounded measures in the sense
that $\mu _{r}\left( D_{r}(G^{\symbol{94}})\right) $ $<$ $\infty $ (and
hence $\mu _{r}(S)$ $<$ $\infty ,$ for all $S$ $\in $ $P\left( D_{r}(G^{%
\symbol{94}})\right) $), and that these measures have more conditions than
any other graph measures introduced in this paper.

\strut

It is also easy to verify that reduced diagram measures are bounded
measures. But, in general, the energy measures and graph groupoid measures
are not bounded measures. However, the energy measures and graph groupoid
measures are locally bounded measures in the sense that $\mu _{G^{\symbol{94}%
}}(F)$ $<$ $\infty $ and $\mu _{\Bbb{G}}(B)$ $<$ $\infty ,$ for all finite
subsets $F$ $\in $ $P\left( \Bbb{F}^{+}(G^{\symbol{94}})\right) $ and $B$ $%
\in $ $P(\Bbb{G}),$ respectively.

\strut

\strut

\strut

\subsection{$\mu _{r}$-Measurable Functions}

\strut

\strut

Let $G$ be the given finite directed graph and $G^{\symbol{94}},$ the
shadowed graph and let $D_{r}(G^{\symbol{94}})$ be the reduced diagram set
of $G^{\symbol{94}}$ and $(G^{\symbol{94}},$ $\mu _{r}),$ the $\mu _{r}$%
-measure space. All simple functions $g$ are defined by

\strut

(2.1) $\ \ \ \ \ \ \ \ \ \ \ \ \ \ \ \ g=\sum_{n=1}^{N}a_{n}1_{S_{n}},$ \
for \ $a_{j}$ $\in $ $\Bbb{C},$

\strut

where $S_{1},$ ..., $S_{N}$ are subsets of $D_{r}(G^{\symbol{94}})$ and where

\strut

\begin{center}
$1_{S_{j}}(w)$ $=$ $\left\{ 
\begin{array}{ll}
1 & \text{if }w\in S_{j} \\ 
0 & \text{otherwise}
\end{array}
\right. $
\end{center}

\strut

are the characteristic functions of $S_{j}$, for all $j$ $=$ $1,$ ..., $N.$
All $\mu _{r}$-measurable functions are approximated by simple functions, in
the same manner of the general measure theory. Recall that all elements in
the reduced diagram set $D_{r}(G^{\symbol{94}})$ are basic elements of the
graph groupoid $\Bbb{G}.$

\strut

\begin{example}
Let $w$ $\in $ $D_{r}(G^{\symbol{94}}).$ Then this element $w$ acts as one
of $\mu _{r}$-measurable functions. Note that $w$ acts on $D_{r}(G^{\symbol{%
94}}),$ as the left multiplication or the right multiplication, in the sense
that $ww^{\prime }$ respectively $w^{\prime }w,$ for all $w^{\prime }$ $\in $
$D_{r}(G^{\symbol{94}}).$ So, we can construct two $\mu _{r}$-measurable
functions $g_{l}^{w}$ and $g_{r}^{w}$ such that

\strut

\ \ \ \ \ \ \ \ \ \ \ \ \ \ \ \ \ \ \ \ \ \ $\ g_{l}^{w}=1_{\delta
^{r}(S_{l}^{w})}$ \ \ and \ \ $g_{r}^{w}=1_{\delta ^{r}(S_{r}^{w})},$

\strut where

\ $\ \ \ \ \ \ \ \ \ \ \ \ \ \ \ S_{l}^{w}=\{w^{\prime }\in \Bbb{F}^{+}(G^{%
\symbol{94}}):ww^{\prime }\in \Bbb{F}^{+}(G^{\symbol{94}})\}$

and

$\ \ \ \ \ \ \ \ \ \ \ \ \ \ \ \ S_{r}^{w}=\{w^{\prime \prime }\in \Bbb{F}%
^{+}(G^{\symbol{94}}):w^{\prime \prime }w\in \Bbb{F}^{+}(G^{\symbol{94}})\}.$

\strut \strut

Therefore, the element $w$ act as a $\mu _{r}$-measurable function $g_{w}$
on $D_{r}(G^{\symbol{94}})$ defined by

\strut

\ \ \ \ \ \ \ \ \ \ \ \ \ \ \ \ \ \ \ \ \ \ $\ \ g_{w}\overset{def}{=}%
g_{l}^{w}+g_{r}^{w}=1_{\delta ^{r}(S_{l}^{w})\cup \delta ^{r}(S_{r}^{w})}.$
\end{example}

\strut \strut \strut

\begin{remark}
For $w$ $\in $ $D(G^{\symbol{94}}),$ we can understand $w$ as $\mu _{\delta
} $-measurable function $g_{w}$ defined by $1_{\delta (S_{l}^{w})\cup \delta
(S_{r}^{w})}.$ For $w$ $\in $ $\Bbb{F}^{+}(G^{\symbol{94}})$ or for $w$ $\in 
$ $\Bbb{G},$ we can determine the $\mu _{G^{\symbol{94}}}$-measurable
function $g_{w}$ as a characteristic function $g_{w}$ $=$ $1_{S_{l}^{w}\cup
S_{r}^{w}}.$ And, for $w$ $\in $ $\Bbb{G},$ the $\mu _{\Bbb{G}}$-measurable
function $g_{w}$ as a characteristic function $g_{w}$ $=$ $1_{S_{l:\Bbb{G}%
}^{w}\cup S_{r:\Bbb{G}}^{w}},$ respectively, where

\strut

$\ \ \ \ \ \ \ \ \ \ \ \ \ \ \ \ \ S_{l:\Bbb{G}}^{w}$ $=$ $S_{l}^{w}$ $\cap $
$\Bbb{G}$ \ \ and \ \ $S_{r:\Bbb{G}}^{w}$ $=$ $S_{r}^{w}$ $\cap $ $\Bbb{G}.$
\end{remark}

\strut \strut

\strut \strut

\strut

\subsection{$\mu _{r}$-Integration}

\strut

\strut

In this section, we will define the integrals of the given $\mu _{r}$%
-measurable functions. First, let $g$ be a simple function given in (2.1).
Then the integral $I_{G}(g)$ of $g$ with respect to $\mu _{r}$ is defined by

\strut

(2.2) \ $\ \ \ \ \ \ \ \ \ I_{G}(g)\overset{denote}{=}\int_{G^{\symbol{94}%
}}g\,\,d\mu _{r}\overset{def}{=}\sum_{n=1}^{N}a_{n}\,\,\mu _{r}(S_{n}).$

\strut

Since $\mu _{r}=d\cup \rho $ $=$ $\rho ,$ where $d$ $=$ $0$ and $\rho $ $=$ $%
\rho $ $\mid _{D_{FP}^{r}(G^{\symbol{94}})},$ the definition (2.2) can be
rewritten as

\strut

\begin{center}
$
\begin{array}{ll}
I_{G}(g) & \,=\sum_{n=1}^{N}a_{n}\mu _{G^{\symbol{94}}}(S_{n}) \\ 
& 
\begin{array}{l}
\\ 
=\sum_{n=1}^{N}a_{n}\left( d(S_{n,V})+\rho (S_{n,FP})\right)
\end{array}
\\ 
& 
\begin{array}{l}
\\ 
=\sum_{n=1}^{N}a_{n}\left( \underset{w\in S_{n,FP}}{\sum }\left| w\right|
\right) ,
\end{array}
\end{array}
$
\end{center}

\strut \strut

where $S_{n,V}$ $=$ $S_{n}$ $\cap $ $V(G^{\symbol{94}})$ and $S_{n,FP}$ $=$ $%
S$ $\cap $ $D_{FP}^{r}(G^{\symbol{94}}),$ for $n$ $=$ $1,$ ..., $N.$

\strut

\begin{remark}
Similar to (2.2), we can define graph integrals for the energy measure $\mu
_{G^{\symbol{94}}},$ the diagram measure $\mu _{\delta }$ and the graph
groupoid measure $\mu _{\Bbb{G}},$ with respect to the $\sigma $-algebras $%
P\left( \Bbb{F}^{+}(G^{\symbol{94}})\right) ,$ $P\left( D(G^{\symbol{94}%
})\right) $ and $P\left( \Bbb{G}\right) ,$ respectively. So,

\strut

$\ \ \ \ \ \ \ \ \ \ \ \ \ \ \int_{G^{\symbol{94}}}$ $g_{1}\,d\mu _{G^{%
\symbol{94}}}=\sum_{n=1}^{N}a_{n}\left( \underset{w\in S_{n,FP}}{\sum }%
\left| w\right| \right) $

and

$\ \ \ \ \ \ \ \ \ \ \ \ \ \ \ \ \int_{\Bbb{G}}\,g\,_{2}\,d\mu _{\Bbb{G}}$ $%
= $ $\sum_{n=1}^{N}a_{n}\left( \underset{w\in S_{n,FP}}{\sum }\left|
w\right| \right) ,$

\strut

whenever $g_{k}$ $=$ $\sum_{n=1}^{N}$ $a_{n}$ $1_{S_{n}^{(k)}}$ are simple
functions, where $a_{n}$ $\in $ $\Bbb{C},$ and $S_{1}^{(1)},$ ..., $%
S_{N}^{(1)}$ $\in $ $P\left( \Bbb{F}^{+}(G^{\symbol{94}})\right) $ and $%
S_{1}^{(2)},$ ..., $S_{N}^{(2)}$ $\in $ $P(\Bbb{G}),$ for $k$ $=$ $1,$ $2,$
respectively.
\end{remark}

\strut

We can easily verify the following proposition.

\strut

\begin{proposition}
Let $g_{1}$ $=\sum_{j=1}^{n}1_{S_{j}}$ and $g_{2}=\sum_{i=1}^{m}1_{T_{i}}$
be simple functions, where $S_{j}$'s and $T_{i}$'s are subsets of $D_{r}(G^{%
\symbol{94}}).$ Suppose that $S_{j}$'s are mutually disjoint and also $T_{i}$%
's are mutually disjoint. If $\cup _{j=1}^{n}$ $S_{j}$ $=$ $\cup _{i=1}^{m}$ 
$T_{i}$ in $D_{r}(G^{\symbol{94}}),$ then $I_{G}(g_{1})$ $=$ $I_{G}(g_{2}).$ 
$\square $
\end{proposition}

$\strut $\strut \strut \strut

\begin{remark}
The above proposition holds true for other cases with respect to other graph
measures. Also, the above proposition holds true when we deal with the
weighted graph measures.
\end{remark}

\strut \strut \strut

Let $1_{S}$ and $1_{T}$ be the characteristic functions, where $S$ $\neq $ $%
T $ in $D_{r}(G^{\symbol{94}}).$ Then we have that $1_{S}$ $\cdot $ $1_{T}$ $%
=$ $1_{S\cap T}.$ Thus we can get the following proposition.

\strut\strut

\begin{proposition}
Let $g_{j}$ $=$ $\sum_{k=1}^{n}$ $a_{j,k}$ $1_{S_{j,k}}$ be a simple
functions, for $j$ $=$ $1,$ $2.$ Then

\strut

$\ \ \ \ I_{G}(g_{1}g_{2})=\sum_{k,i=1}^{n}\left( a_{1,k}a_{2,i}\right)
\cdot \left( \underset{w\in S_{(1,k),FP}\cap S_{(2,i),FP}}{\sum }\left|
w\right| \right) .$
\end{proposition}

\strut

\begin{proof}
Observe that

\strut

$\ \ \ g_{1}g_{2}=\left( \sum_{k=1}^{n}a_{1,k}1_{S_{1,k}}\right) \left(
\sum_{i=1}^{n}a_{2,i}1_{S_{2,i}}\right) $

\strut

$\ \ \ \ \ \ \ \ \ \ \ =\sum_{k,i=1}^{n}a_{1,k}a_{2,i}\left(
1_{S_{1,k}}\cdot1_{S_{2,i}}\right)
=\sum_{k,i=1}^{n}a_{1,k}a_{2,i}1_{S_{1,k}\cap S_{2,i}}$

\strut

So, we have that

\strut

(2.3) $\ \ \ \ \ \ \ I_{G}(g_{1}g_{2})=\sum_{k,i=1}^{n}a_{1,k}a_{2,i}\mu
_{r}\left( S_{1,k}\cap S_{2,i}\right) .$

\strut

Consider

\strut

$\ \ \mu _{r}\left( S_{1,k}\cap S_{2,i}\right) =d\left( (S_{1,k}\cap
S_{2,i})_{V}\right) +\rho \left( (S_{1,k}\cap S_{2,i})_{FP}\right) $

\strut

$\ \ \ \ \ \ \ \ \ \ \ \ \ \ \ \ \ \ \ \ \ \ \ \ \ \ \ =\rho \left(
S_{(1,k),FP}\cap S_{(2,k),FP}\right) $

\strut

$\ \ \ \ \ \ \ \ \ \ \ \ \ \ \ \ \ \ \ \ \ \ \ \ \ \ \ =\underset{w\in
S_{(1,k),FP}\cap S_{(2,i),FP})}{\sum }\left| w\right| .$

$\strut$
\end{proof}

\strut

The formula (2.3) is satisfied for all other graph measures when we replace $%
\mu _{r}$ to $\mu _{G^{\symbol{94}}},$ $\mu _{\Bbb{G}}$ and $\mu _{\delta }$%
. Suppose that $g_{1}$ and $g_{2}$ are given as in the previous proposition
and assume that the families $\{S_{1,k}\}_{k=1}^{n}$ and $%
\{S_{2,i}\}_{i=1}^{n}$ are disjoint in $D_{r}(G^{\symbol{94}})$. i.e.,

$\strut$

\begin{center}
$\left( \cup_{k=1}^{n}S_{1,k}\right) \cap\left( \cup_{i=1}^{n}S_{2,i}\right)
=\varnothing.$
\end{center}

\strut

Then $I(g_{1}g_{2})=0.$

\strut \strut \strut

If $g$ is a $\mu _{r}$-measurable function, then the support of $g$ is
denoted by $D_{r}(G^{\symbol{94}}$ $:$ $g)$\strut , as a subset of $D_{r}(G^{%
\symbol{94}}).$ i.e.,

\strut

\begin{center}
$D_{r}(G^{\symbol{94}}:g)\overset{def}{=}\{w\in D_{r}(G^{\symbol{94}%
}):g(w)\neq 0\}.$
\end{center}

\strut \strut

Let's observe more $\mu _{r}$-measurable functions. First, consider the
monomial $g_{1}(x)$ $\overset{def}{=}$ $g_{x}$ $=$ $g_{\delta _{x}^{r}},$
for all $x$ $\in $ $D_{r}(G^{\symbol{94}}),$ where $g_{x}$ $=$ $1_{\delta
^{r}(S_{l}^{x})\cup \delta ^{r}(S_{r}^{x})},$ for all $x$ $\in $ $D_{r}(G^{%
\symbol{94}}).$ Then the support $D_{r}(G^{\symbol{94}}$ $:$ $g_{1})$ is
same as the reduced diagram set $D_{r}(G^{\symbol{94}}),$ since $g_{x}$ is
well-defined on $D_{r}(G^{\symbol{94}}),$ for all $x$ $\in $ $D_{r}(G^{%
\symbol{94}}).$ Thus we have that

\strut

(2.4) $\ \ \ \ \ \ \ \ \ \ \ I_{G}(g_{1})=\underset{w\in D_{r}(G^{\symbol{94}%
})}{\sum }\mu _{r}\left( \delta ^{r}(S_{l}^{w})\cup \delta
^{r}(S_{r}^{w}))\right) $.

\strut

Notice that if $G$ is a finite directed graph, then $D_{r}(G^{\symbol{94}})$
is a finite set. i.e., $\left| D_{r}(G^{\symbol{94}})\right| $ $<$ $\infty .$
Therefore, we can have that;

\strut

\begin{proposition}
If $g_{1}$ is a $\mu _{r}$-measurable function defined by $g_{1}(w)$ $%
\overset{def}{=}$ $g_{w},$ for all $w$ $\in $ $D_{r}(G^{\symbol{94}}),$ then 
$I_{G}(g_{1})$ $<$ $\infty .$ $\square $
\end{proposition}

\strut

\begin{remark}
Suppose $g_{1}$ is similarly defined as a $\mu _{\delta }$-measurable
function, where $\mu _{\delta }$ is the diagram measure. Then we can get the
same result as in the previous proposition. i.e., $\int_{G^{\symbol{94}}}$ $%
g_{1}$ $d\,\mu _{\delta }$ $<$ $\infty .$ However, if we define $g_{1}$
similarly with respect to the energy measure $\mu _{G^{\symbol{94}}}$ or
with respect to the graph groupoid measure $\mu _{\Bbb{G}},$ then we cannot
guarantee the boundedness of $\int_{G^{\symbol{94}}}$ $g_{1}$ $d\,\mu _{G^{%
\symbol{94}}}$ and $\int_{G^{\symbol{94}}}$ $g_{1}$ $d\,\mu _{\Bbb{G}}.$ For
instance, if the graph $G$ is a one-vertex-2-loop-edge graph $G_{2},$ then
the integrals $\int_{G_{2}^{\symbol{94}}}$ $g_{1}$ $d\mu _{G_{2}^{\symbol{94}%
}}$ $=$ $\infty $ $=$ $\int_{G_{2}^{\symbol{94}}}$ $g_{1}$ $d$ $\mu _{\Bbb{G}%
_{2}^{\symbol{94}}},$ since the supports

\strut

$\ \ \ \ \ \ \ \Bbb{F}^{+}(G_{2}^{\symbol{94}}:g_{1})=\Bbb{F}^{+}(G_{2}^{%
\symbol{94}})$ \ and \ $\Bbb{G}_{2}(G^{\symbol{94}}:g_{1})=\Bbb{G}_{2},$

\strut

respectively and

\strut

$\ \ \ \ \ \ \ \ \ \ \ \ \ \ \ \ \ \ \ \left| \Bbb{F}^{+}(G_{2}^{\symbol{94}%
})\right| =\infty =\left| \Bbb{G}_{2}\right| ,$

\strut
\end{remark}

\strut

\ Now, let $g_{2}(x)\overset{def}{=}g_{\delta _{x^{2}}^{r}}$ on $D_{r}(G^{%
\symbol{94}}).$ This monomial $g_{2}$ has its support

\strut

\begin{center}
$D_{r}(G^{\symbol{94}}:g_{2})=\left\{ w\in D_{r}(G^{\symbol{94}}):\left| 
\begin{array}{c}
w\in V(G)\text{ \ \ \ \ \ \ \ \ \ \ \ or} \\ 
w\text{ is a loop finite path}
\end{array}
\right. \right\} .$
\end{center}

\strut

The support $D_{r}(G^{\symbol{94}}:g_{2})$ is determined as above, since $%
g_{2}(v)$ $=$ $g_{v^{2}}$ $=$ $g_{v},$ for all $v$ $\in $ $V(G^{\symbol{94}%
}),$ and $g_{2}(l)$ $=$ $g_{\delta _{l^{2}}^{r}}$ $=$ $g_{l},$ for all
(basic) loop finite paths $l$ in $D_{r}(G).$ Assume now that $w$ $=$ $v_{1}$ 
$w$ $v_{2}$ is a non-loop finite path in $D_{r}(G^{\symbol{94}})$ with $%
v_{1} $ $\neq $ $v_{2}$ in $V(G^{\symbol{94}}).$ Then $w^{2}$ $=$ $(v_{1}$ $%
w $ $v_{2})$ $(v_{1}$ $w$ $v_{2})$ $=$ $\emptyset ,$ and hence $\delta
_{w^{2}} $ $=$ $\delta _{\emptyset }$ $=$ $\emptyset ,$ in $D_{r}(G^{\symbol{%
94}})$\strut $.$ So, in this case, $w$ $\notin $ $D_{r}(G^{\symbol{94}}$ $:$ 
$g_{2}).$ We can get that

\strut

(2.5) \ 

\begin{center}
$
\begin{array}{ll}
I_{G}(g_{2}) & =\underset{w\in D_{r}(G^{\symbol{94}}:g_{2})}{\sum }\mu
_{r}\left( \delta ^{r}(S_{l}^{w})\cup \delta ^{r}(S_{r}^{w})\right) \\ 
&  \\ 
& =\,\underset{w\text{ : loop in }D_{r}(G^{\symbol{94}})}{\sum }\mu
_{r}\left( \delta ^{r}(S_{l}^{w})\cup \delta ^{r}(S_{r}^{w})\right) .
\end{array}
$
\end{center}

\strut

In general, we can get that;

\strut

\begin{proposition}
Let $g_{n}(x)\overset{def}{=}g_{\delta _{x^{n}}^{r}}$ be the monomial, for
all $x$ $\in $ $D_{r}(G^{\symbol{94}})$. Then $I_{G}\left( g_{n}\right) $ $=$
$I_{G}\left( g_{2}\right) ,$ for all $n$ $\in $ $\mathbb{N}$ $\setminus $ $%
\{1\},$ where $I_{G}(g_{2})$ satisfies (2.5).
\end{proposition}

\strut

\begin{proof}
It suffices to show that the support $D_{r}(G^{\symbol{94}}:g_{n})$ of $%
g_{n} $ and the support $D_{r}(G:g_{2})$ of $g_{2}$ coincide, for all $n$ $%
\in $ $\mathbb{N}$ $\setminus $ $\{1\}.$ It is easy to check that if $w$ is
in $D_{r}(G^{\symbol{94}}),$ then $w^{n}$ exists in $\Bbb{F}^{+}(G^{\symbol{%
94}}) $ if and only if either $w$ is a vertex or $w$ is a loop finite path,
for $n$ $\in $ $\Bbb{N}$ $\setminus $ $\{1\}.$ So, the support of $g_{n}$ is

\strut

$\ \ \ \ \ \ \ \ \ \ \ D_{r}(G^{\symbol{94}}:g_{n})=V(G)\cup loop_{r}(G),$
for $n$ $\in $ $\Bbb{N}$ $\setminus $ $\{1\}.$

\strut

where $loop_{r}(G)\overset{def}{=}\{l\in D_{r}(G):l$ is a loop finite path$%
\}.$ Therefore, the support $D_{r}(G^{\symbol{94}}$ $:$ $g_{n})$ is same as $%
D_{r}(G^{\symbol{94}}$ $:$ $g_{2})$ whenever $n$ $\neq $ $1.$ Therefore,

\strut

$\ \ \ \ \ \ \ \ \ \ \ \ \ 
\begin{array}{ll}
I_{G}\left( g_{n}\right) & =\underset{w\in D_{r}(G^{\symbol{94}}:g_{2})}{\sum%
}I_{G}(g_{w}) \\ 
&  \\ 
& =\underset{w\in D_{r}(G^{\symbol{94}}:g_{2})}{\sum}\mu_{G^{\symbol{94}%
}}\left( \delta^{r}(S_{l}^{w})\cup\delta^{r}(S_{r}^{w})\right) .
\end{array}
$

$\strut\strut$
\end{proof}

\strut

\begin{remark}
The monomial $g_{n}$ on $\Bbb{F}^{+}(G^{\symbol{94}})$ and on $\Bbb{G}$ are
defined by $g_{n}(w)$ $=$ $g_{w^{n}},$ for all $w$ $\in $ $\Bbb{F}^{+}(G^{%
\symbol{94}}),$ respectively, for all $w$ $\in $ $\Bbb{G}.$ The monomial $%
g_{n}$ on $D(G^{\symbol{94}})$ is defined by $g_{n}(w)$ $=$ $g_{\delta
_{w^{n}}}$, for all $w$ $\in $ $D(G^{\symbol{94}}).$ When we deal with the
energy measure $\mu _{G^{\symbol{94}}}$ and the graph groupoid measure $\mu
_{\Bbb{G}},$ the above proposition holds, similarly, but the meaning is of
course different. For example, on $D_{r}(G^{\symbol{94}})$ or on $D(G^{%
\symbol{94}}),$ there exists $t$ $<$ $\infty $ in $\Bbb{N}$ $\cup $ $\{0\}$
such that the cardinality of the set of all loops in $D_{r}(G^{\symbol{94}})$
or in $D(G^{\symbol{94}})$ is $t,$ since $D_{r}(G^{\symbol{94}})$
respectively $D(G^{\symbol{94}})$ is a finite set, by the finiteness of $G.$
Therefore,

\strut

$\ \ \ I_{G}(g_{n})$ $=$ $I_{G}(g_{2})$ $<$ $\infty $ and $\int_{G^{\symbol{%
94}}}$ $g_{n}$ $d\mu _{\delta }$ $=$ $\int_{G^{\symbol{94}}}g_{n}\,d\mu
_{\delta }$ $<$ $\infty .$

\strut

But on the graph groupoid $\Bbb{G},$ the set of all loop finite paths will
be $\{l_{i}^{k}$ $:$ $i$ $=$ $1,$ ..., $t,$ $k$ $\in $ $\Bbb{Z}$ $\setminus $
$\{0\}\},$ whenever $loop_{r}$ $(G^{\symbol{94}})$ $=$ $\{l_{1},$ ..., $%
l_{t}\}.$ Thus the graph integral

\strut

$\ \ \ \ \ \ \ \ \ \ \ \ \ \ \ \ \ \int_{G^{\symbol{94}}}$ $g_{n}$ $d\mu _{%
\Bbb{G}}$ $=$ $\infty ,$ for all $n$ $\in $ $\Bbb{N}$,

\strut

whenever $t$ $\geq $ $1.$ Similarly, $\int_{G^{\symbol{94}}}$ $g_{n}$ $d\mu
_{G^{\symbol{94}}}$ $=$ $\infty .$
\end{remark}

\strut

In the following theorem, we compute the integral of polynomials.\strut

\strut \strut

\begin{theorem}
Let $g_{p}$ $\overset{def}{=}$ $\sum_{n=0}^{N}a_{n}g_{n}$ be a polynomial
with $g_{0}$ $\overset{def}{=}$ $1,$ for $a_{0},$ ..., $a_{N}$ $\in $ $%
\mathbb{C}.$ Then

\strut

\strut (2.6)

$\ \ \ \ 
\begin{array}{ll}
I_{G}\left( g_{p}\right) = & a_{0}\mu _{r}\left( D_{r}(G^{\symbol{94}%
})\right) +a_{1}\left( \underset{w\in D_{r}(G^{\symbol{94}})}{\sum }\mu
_{r}\left( \delta ^{r}(S_{l}^{w})\cup \delta ^{r}(S_{r}^{w})\right) \right)
\\ 
&  \\ 
& +\left( \underset{w\in V(G)\cup loop_{r}(G)}{\sum }\mu _{r}\left( \delta
^{r}(S_{l}^{w})\cup \delta ^{r}(S_{r}^{w})\right) \right) \left(
\sum_{k=2}^{N}a_{k}\right) ,
\end{array}
$

\strut

where $loop_{r}(G)=\{l\in D_{r}(G^{\symbol{94}}):l$ is a loop finite path$%
\}. $
\end{theorem}

\strut

\begin{proof}
Let $g_{p}$ be the given polynomial on $\mathbb{F}^{+}(G^{\symbol{94}}).$
Then

\strut

$\ \ \ I_{G}\left( g_{p}\right) =\sum_{n=0}^{N}a_{n}I_{G}\left( g_{n}\right) 
$

\strut

$\ \ \ \ \ \ \ =a_{0}I_{G}(1)+a_{1}I_{G}\left( g_{1}\right)
+\sum_{k=2}^{N}a_{k}I_{G}\left( g_{k}\right) $

\strut

$\ \ \ \ \ \ \ =a_{0}\mu _{r}\left( D_{r}(G^{\symbol{94}})\right)
+a_{1}\left( \underset{w\in D_{r}(G^{\symbol{94}})}{\sum }\mu _{r}(\delta
^{r}(S_{l}^{w})\cup \delta ^{r}(S_{r}^{w}))\right) $

\strut

$\ \ \ \ \ \ \ \ \ \ \ \ \ \ \ \ \ \ \ \ \ \ \ \ \ \ \ \ \ \ \ \ \ \ \ \ \ \
\ +\sum _{k=2}^{N}a_{k}I_{G}\left( g_{2}\right) $

\strut

since the constant function $1$ has its support, $D_{r}(G^{\symbol{94}}),$
and since $I_{G}\left( g_{n}\right) $ $=$ $I_{G}\left( g_{2}\right) ,$ for
all $n$ $\geq $ $2$

\strut

$\ \ \ \ \ \ \ =a_{0}\mu _{r}\left( D_{r}(G)\right) +a_{1}\left( \underset{%
w\in D_{r}(G^{\symbol{94}})}{\sum }\mu _{r}(\delta ^{r}(S_{l}^{w})\cup
\delta ^{r}(S_{r}^{w}))\right) $

\strut

$\ \ \ \ \ \ \ \ \ \ \ \ \ \ \ \ \ \ \ \ \ \ \ \ \ \ \ \ \ \ \
+\sum_{k=2}^{N}a_{k}\left( \underset{w\in V(G)\cup loop(G^{\symbol{94}})}{%
\sum }\mu _{r}\left( \delta ^{r}(S_{l}^{w})\cup \delta
^{r}(S_{r}^{w})\right) \right) .$

\strut
\end{proof}

\strut \strut \strut

\begin{remark}
By the previous theorem and by the previous remark, we verify that $\int_{%
\Bbb{G}}$ $g_{p}$ $d\mu _{\Bbb{G}}$ and $\int_{\Bbb{F}^{+}(G^{\symbol{94}})}$
$g_{p}$ $d\mu _{G^{\symbol{94}}}$ have the similar formuli like (2.6).
However, these integrals are $\infty ,$ in general.
\end{remark}

\strut \strut

\begin{corollary}
Let $w$ $=v_{1}wv_{2}$ be a finite path in $D_{r}(G^{\symbol{94}})$ with $%
v_{1}$ $\neq $ $v_{2}$ in $V(G)$, and let $g$ $=$ $\sum_{n=0}^{N}$ $%
a_{n}g_{\delta _{w^{n}}^{r}},$ with $g_{w^{0}}$ $\overset{def}{=}$ $1,$
where $a_{0},$ ..., $a_{N}$ $\in $ $\mathbb{R}.$ Then

\strut

$\ \ \ \ \ \ \ \ I_{G}\left( g\right) =a_{0}\mu _{r}\left( D_{r}(G^{\symbol{%
94}})\right) +a_{1}\mu _{r}\left( \delta ^{r}(S_{l}^{w})\cup \delta
^{r}(S_{r}^{w})\right) .$
\end{corollary}

\strut

\begin{proof}
Since $w$ is a non-loop finite path, $w^{k}$ $=$ $\emptyset,$ for all $k$ $%
\in$ $\mathbb{N}$ $\setminus$ $\{1\}.$ Therefore, $S_{l}^{w^{k}}$ $\cup$ $%
S_{r}^{w^{k}}$ $=$ $\varnothing,$ for all $k$ $=$ $2,$ $3,$ ..., $N.$ This
shows that $I_{G}\left( g_{w}^{k}\right) $ $=$ $0,$ for all $k$ $=$ $2,$ $3,$
..., $N.$ So, $I_{G}\left( g\right) $ $=$ $a_{0}$ $I_{G}\left( 1\right) $ $+$
$a_{1}$ $I_{G}\left( g_{w}\right) .$
\end{proof}

\strut

\begin{corollary}
Let $w$ $=$ $vwv$ be a loop finite path in $D_{r}(G^{\symbol{94}}),$ with $v$
$\in $ $V(G^{\symbol{94}}),$ and let $g$ $=$ $\sum_{n=0}^{N}$ $%
a_{n}g_{\delta _{w^{n}}^{r}}$, with $g_{w^{0}}$ $=$ $1,$ where $a_{0},$ ..., 
$a_{n}$ $\in $ $\mathbb{R}.$ Then

\strut

$\ \ \ \ \ \ I_{G}\left( g\right) =a_{0}\mu _{r}\left( D_{r}(G)\right)
+\left( \delta ^{r}(S_{l}^{w})\cup \delta ^{r}(S_{r}^{w})\right) \left(
\sum_{k=1}^{N}a_{k}\mu _{r}\right) .$

$\square $
\end{corollary}

\strut\strut\strut\strut

Now, we will consider the map $g_{-1}$ defined by $g_{-1}(x)$ $\overset{def}{%
=}$ $g_{x^{-1}}$ $=$ $g_{\delta _{x^{-1}}^{r}}.$ On $D_{r}(G^{\symbol{94}}),$
the map $g_{-1}$ is well-determined because, for any $w$ $\in $ $D_{r}(G^{%
\symbol{94}}),$ there always exists $w^{-1}$ in $D_{r}(G^{\symbol{94}}).$
Moreover, the support $D_{r}(G^{\symbol{94}}$ $:$ $g_{w})$ of $g_{w}$ and
the support $D_{r}(G^{\symbol{94}}$ $:$ $g_{w^{-1}})$ of $g_{w^{-1}}$ are
same. In fact, $D_{r}(G^{\symbol{94}}$ $:$ $g_{w^{-1}})$ $=$ $D_{r}(G^{%
\symbol{94}}$ $:$ $g_{w})^{-1},$ where $X^{-1}$ means the set $\{x^{-1}$ $:$ 
$x$ $\in $ $X\}.$ Indeed,

\strut

\strut (2.7)

\begin{center}
$
\begin{array}{ll}
D_{r}(G^{\symbol{94}}:g_{w}) & =S_{l}^{w}\cup S_{r}^{w}=D_{r}(G^{\symbol{94}%
}:g_{w})^{-1} \\ 
&  \\ 
& =S_{r}^{w^{-1}}\cup S_{l}^{w^{-1}}=D_{r}(G^{\symbol{94}}:g_{w^{-1}}).
\end{array}
$
\end{center}

\strut \strut

Hence, we can get that

\strut

(2.8)\ \ \ \ \ \ \ \ \ \ \ $D_{r}(G^{\symbol{94}}:g_{-1})=D_{r}(G^{\symbol{94%
}})=D_{r}(G^{\symbol{94}}:g_{1}),$

\strut

since $g_{1}$ has its support $D_{r}(G^{\symbol{94}})$.

\strut

\begin{proposition}
Let $g_{-1}$ be given as above. Then $g_{-1}$ is $\mu _{r}$-measurable and $%
I_{G}\left( g_{-1}\right) $ $=$ $I_{G}\left( g_{1}\right) .$ $\square $
\end{proposition}

\strut \strut \strut

\begin{remark}
Define the energy measurable function $g_{-1}$ by $g_{-1}(w)$ $\overset{def}{%
=}$ $g_{w^{-1}},$ for all $w$ $\in $ $\Bbb{F}^{+}(G^{\symbol{94}}),$ and
define the graph groupoid measurable function $g_{-1}$ by $g_{-1}(w)$ $%
\overset{def}{=}$ $g_{w^{-1}},$ for all $w$ $\in $ $\Bbb{G},$ and define the
diagram measurable function $g_{-1}$ by $g_{-1}(w)$ $=$ $g_{\delta
_{w^{-1}}},$ for all $w$ $\in $ $D(G^{\symbol{94}}).$ Then we can have that $%
\Bbb{F}^{+}(G^{\symbol{94}}$ $:$ $g_{-1})$ $=$ $\Bbb{F}^{+}$ $(G^{\symbol{94}%
}$ $:$ $g_{1})$ and $\Bbb{G}(G^{\symbol{94}}$ $:$ $g_{-1})$ $=$ $\Bbb{G}(G^{%
\symbol{94}}$ $:$ $g_{1})$ and $D(G^{\symbol{94}}$ $:$ $g_{-1})$ $=$ $D(G^{%
\symbol{94}}$ $:$ $g_{1}),$ similar to (2.7). Therefore, like the previous
proposition,

\strut

$\ \ \ \ \ \ \int_{G^{\symbol{94}}}g_{-1}\,d\mu _{G^{\symbol{94}}}$ $=$ $%
\int_{G^{\symbol{94}}}g_{1}\,d\mu _{G^{\symbol{94}}},$ \ \ $\int_{G^{\symbol{%
94}}}g_{-1}d\mu _{\Bbb{G}}$ $=$ $\int_{G^{\symbol{94}}}\,g_{1}\,d\mu _{\Bbb{G%
}}$

and

$\ \ \ \ \ \ \ \ \ \ \ \ \ \ \ \ \ \ \ \ \int_{G^{\symbol{94}}}$ $g_{-1}$ $%
d\mu _{\delta }=\int_{G^{\symbol{94}}}g_{1}\,d\mu _{\delta }.$
\end{remark}

\strut \strut

Similar to the previous proposition, we can conclude that;

\strut

\begin{theorem}
Define $g_{-n}(x)\overset{def}{=}g_{\delta _{x^{-n}}^{r}},$ for $n\in %
\mathbb{N}$, for $x$ $\in $ $D_{r}(G^{\symbol{94}}).$ Then $I_{G}\left(
g_{-n}\right) $ $=$ $I_{G}\left( g_{n}\right) $.
\end{theorem}

\strut

\begin{proof}
Observe that if $n$ $\geq$ $2,$ then

\strut

$\ \ D_{r}(G^{\symbol{94}}:g_{-n})=D_{r}(G^{\symbol{94}}:g_{-2})=V(G^{%
\symbol{94}})\cup loop_{r}(G^{\symbol{94}})=D_{r}(G^{\symbol{94}}:g_{2}).$

\strut

So, $I_{G}\left( g_{-n}\right) =I_{G}\left( g_{-2}\right) =I_{G}\left(
g_{2}\right) .$ By the previous proposition, $I_{G}\left( g_{-1}\right) $ $=$
$I_{G}\left( g_{1}\right) .$
\end{proof}

\strut

So, we can consider the trigonometric polynomials on $D_{r}(G^{\symbol{94}%
}). $

\strut

\begin{corollary}
Let $g_{k}$ and $g_{-k}$ be given as above, for all $k$ $\in $ $\mathbb{N},$
and let $g$ $=$ $\sum_{n=-M}^{N}$ $a_{n}g_{n}$ be a trigonometric polynomial
with $g_{0}$ $\overset{def}{=}$ $1,$ where $N,$ $M$ $\in $ $\mathbb{N}.$ Then

\strut \strut

(2.9) $\ I_{G}\left( g\right) =a_{0}\mu _{r}\left( D_{r}(G)\right) +\left(
a_{1}+a_{-1}\right) I_{G}(g_{1})$

\strut

$\ \ \ \ \ \ \ \ \ \ \ \ \ \ \ \ \ \ \ \ \ \ \ \ \ \ \ \ \ \ \ \
+I_{G}\left( g_{2}\right) \cdot \left(
\sum_{n=-M}^{-2}a_{n}+\sum_{k=2}^{N}a_{k}\right) .$
\end{corollary}

\strut

\begin{proof}
Let $g$ be the given trigonometric polynomial on $D_{r}(G^{\symbol{94}}).$
Then

\strut

$I_{G}\left( g\right) =I_{G}\left( \sum_{n=-N}^{N}a_{n}g_{n}\right)
=\sum_{n=-N}^{N}a_{n}I_{G}\left( g_{n}\right) $

\strut

\ $=\sum_{n=-M}^{-2}a_{n}I_{G}\left( g_{n}\right)
+a_{-1}I_{G}(g_{-1})+a_{0}I_{G}\left( g_{0}\right)
+a_{1}I_{G}(g_{1})+\sum_{k=2}^{N}a_{k}I_{G}\left( g_{k}\right) $

\strut

\ $=\sum_{n=-M}^{-2}a_{n}I_{G}\left( g_{n}\right) +a_{0}\mu_{G^{\symbol{94}%
}}\left( D_{r}(G)\right) ++\sum_{k=1}^{N}a_{k}I_{G}\left( g_{k}\right) $

\strut

by the previous proposition

\strut

$=a_{0}\mu _{r}\left( D_{r}(G)\right)
+(a_{-1}+a_{1})I_{G}(g_{1})+\sum_{n=-M}^{-1}a_{n}I_{G}\left( g_{2}\right)
+\sum_{k=1}^{N}a_{k}I_{G}\left( g_{2}\right) $

\strut

by the fact that $I_{G}\left( g_{n}\right) =I_{G}\left( g_{2}\right) ,$ for
all $n$ $\in$ $\mathbb{N}$ $\setminus$ $\{1\}.$
\end{proof}

\strut\strut\strut

\strut\strut

\strut\strut

\subsection{Examples}

\strut

\strut

In this section, we will consider certain finite directed graphs and
corresponding graph measures and graph integrals. Let $G_{\Lambda }$ be a
tree with

\strut

\begin{center}
$V(G_{\Lambda })=\{v_{1},v_{2},v_{3}\}$ \ and \ $E(G_{\Lambda })=\{e_{1}$ $=$
$v_{1}e_{1}v_{2},$ $e_{2}=v_{1}e_{2}v_{3}\}.$
\end{center}

\strut

Let $G_{\Delta }$ be a one-flow circulant graph with

\strut

\begin{center}
$V(G_{\Delta })=\{v_{1},v_{2},v_{3}\}$ \ and \ $E(G_{\Delta })=\left\{ 
\begin{array}{c}
e_{1}=v_{1}e_{1}v_{2}, \\ 
e_{2}=v_{2}e_{2}v_{3}, \\ 
e_{3}=v_{3}e_{3}v_{1}
\end{array}
\right\} .$
\end{center}

\strut \strut \strut \strut \strut \strut

\begin{example}
Consider $G_{\Lambda }.$ We have that

\strut

$\ \ \ \ \ \ \ \ \ \ \ 
\begin{array}{ll}
I_{G_{\Lambda }}\left( g_{v_{1}}\right) & \,=\mu _{r}\left(
\{v_{1},e_{1},e_{2},e_{1}^{-1},e_{2}^{-1}\}\right) \\ 
& \,=d(\{v_{1}\})+\rho \left( \{e_{1},e_{2},e_{1}^{-1},e_{2}^{-1}\}\right)
\\ 
& 
\begin{array}{l}
=0+\rho \left( e_{1}\right) +\rho (e_{2})+\rho \left( e_{1}^{-1}\right)
+\rho \left( e_{2}^{-1}\right) \\ 
=4.
\end{array}
\end{array}
$

\strut

$\ \ \ \ \ \ \ \ 
\begin{array}{ll}
I_{G_{\Lambda }}\left( g_{v_{2}}\right) & =\mu
_{r}(\{v_{2},e_{1},e_{1}^{-1}\})=d(\{v_{2}\})+\rho (\{e_{1},e_{1}^{-1}\}) \\ 
& =0+2=2.
\end{array}
$

\strut

Similarly, $I_{G_{\Lambda }}\left( g_{v_{3}}\right) =2.$

\strut

$\ \ \ \ \ \ \ \ \ \ \ \ \ \ \ 
\begin{array}{ll}
I_{G_{\Lambda }}\left( g_{e_{1}}\right) & \,=\mu _{G_{\Lambda }^{\symbol{94}%
}}\left( \{v_{1},v_{2},e_{1}^{-1},e_{1}^{-1}e_{2},\,e_{2}^{-1}\}\right) \\ 
& 
\begin{array}{l}
=d\left( \{v_{1},v_{2}\}\right) +\rho (e_{1}^{-1},e_{1}^{-1}e_{2},e_{2}^{-1})
\\ 
=0+1+2+1=4=I_{G_{\Lambda }}\left( g_{e_{1}^{-1}}\right) .
\end{array}
\end{array}
$

\strut

Similarly, $I_{G_{\Lambda }}\left( g_{e_{2}}\right) =4=I_{G_{\Lambda
}}\left( g_{e_{2}^{-1}}\right) .$ We have that

\strut

$\ \ \ \ \ \ \ \ D_{r}(G_{\Lambda })=\{v_{1},v_{2},v_{3}\}\cup
\{e_{1},e_{2},e_{1}^{-1},e_{2}^{-1},e_{1}^{-1}e_{2},e_{2}^{-1}e_{1},\}.$

\strut

Thus

\strut

$\ \ \ \ 
\begin{array}{ll}
I_{G_{\Lambda }}\left( g_{1}\right) & =\underset{w\in D_{r}(G_{\Lambda }^{%
\symbol{94}})}{\sum }I_{G_{\Lambda }}\left( g_{w}\right) \\ 
& 
\begin{array}{l}
=\sum_{j=1}^{3}I_{G_{\Lambda }}\left( g_{v_{j}}\right)
+\sum_{k=1}^{2}I_{G_{\Lambda }}\left( g_{e_{k}}\right) \\ 
\text{ \ \ }+\sum_{i=1}^{2}I_{G_{\Lambda }}\left( g_{e_{i}^{-1}}\right)
+I_{G_{\Lambda }}(g_{e_{1}^{-1}e_{2}})+I_{G_{\Lambda }}(g_{e_{2}^{-1}e_{1}})
\end{array}
\\ 
& =\left( 4+2+2\right) +4+4+\left( 0+0\right) =16.
\end{array}
$

$\strut \strut $

We also have that

\strut

$\ \ \ \ \ \ \ \ \ \ \ \ I_{G_{\Lambda }}\left( g_{n}\right)
=\sum_{j=1}^{3}I_{G_{\Lambda }}\left( g_{v_{j}}\right) =4+2+2=8,$

\strut

for all $n$ $\in $ $\mathbb{N}$ $\setminus $ $\{1\},$ because $w^{k}$ $=$ $%
\emptyset ,$ for all $w$ $\in $ $D_{FP}^{r}(G_{\Lambda }^{\symbol{94}})$ and 
$k$ $\in $ $\mathbb{N}$ $\setminus $ $\{1\},$ and $D_{r}(G_{\Lambda }^{%
\symbol{94}}$ $:$ $g_{k})$ $=$ $V(G_{\Lambda }^{\symbol{94}})$ $\cup $ $%
loop_{r}(G_{\Lambda }^{\symbol{94}})$ $=$ $V(G_{\Lambda }^{\symbol{94}}).$
Therefore,

\strut

$\ \ \ \ \ \ \ \ \ \ \ \ \ \ \ \ \ I_{G_{\Lambda }}\left( g_{p}\right)
=1+16+8N=8N+17,$

\strut

where $g_{p}=\sum_{n=0}^{N}g_{n}$, with $g_{0}$ $=$ $1,$ for all $N$ $\in $ $%
\mathbb{N}.$ Also,

\strut

$\ \ \ \ \ \ \ \ I_{G_{\Lambda }}\left( g_{t}\right)
=8M+16+1+16+8N=8(N+M)+33,$

\strut

where $g_{t}=\sum_{n=-M}^{N}g_{n},$ with $g_{0}$ $=$ $1,$ for all $N,$ $M$ $%
\in $ $\Bbb{N}.$
\end{example}

\strut \strut \strut \strut

\begin{example}
Consider the graph $G_{\Delta }.$ Then we have that

\strut

$\ \ \ 
\begin{array}{ll}
I_{G_{\Delta }}\left( g_{v_{1}}\right) & =\mu _{r}\left( \left\{ 
\begin{array}{c}
v_{1},\,e_{1}^{\pm 1},e_{3}^{\pm
1},e_{1}e_{2},e_{2}^{-1}e_{1}^{-1},\,e_{2}e_{3}, \\ 
\,\,\,\,\,e_{3}^{-1}e_{2}^{-1},e_{1}e_{2}e_{3},e_{3}^{-1}e_{2}^{-1}e_{1}^{-1}
\end{array}
\right\} \right) \\ 
& =18
\end{array}
$

\strut

Similarly, $I_{G_{\Delta }}\left( g_{v_{2}}\right) =18=I_{G_{\Delta }}\left(
g_{v_{3}}\right) .$ Also we can have that

\strut

\ \ \ \ \ \ \ \ \ \ \ \ $I_{G_{\Delta }}\left( g_{e_{j}}\right)
=4=I_{G_{\Delta }}\left( g_{e_{j}^{-1}}\right) ,$ for all $j$ $=$ $1,$ $2,$ $%
3.$

\strut

Indeed, if $j$ $=$ $1,$ then $I_{G_{\Delta }}\left( g_{e_{1}}\right) =\mu
_{G_{\Delta }^{\symbol{94}}}\left( \delta ^{r}(S_{l}^{w})\cup \delta
^{r}(S_{r}^{w})\right) ,$ where

\strut

$\ \ \ \ \ \ \delta ^{r}\left( S_{l}^{e_{1}}\right)
=\{v_{2},e_{2},e_{2}e_{3}\}$ \ and \ $\delta ^{r}\left( \text{ }%
S_{r}^{e_{1}}\right) =\{v_{1},e_{3},e_{2}e_{3}\},$

so,

$\ \ \ \ \ \ \ \ \ \ \ \ \ \ \ \ \ \ \ \ S_{l}^{e_{1}}\cup
S_{r}^{e_{1}}=\{v_{1},v_{2},e_{2},e_{3},\,e_{2}e_{3}\}.$

\strut Thus

\strut

$\ \ \ I_{G_{\Delta }}\left( g_{e_{1}}\right) =\mu _{G_{\Delta }^{\symbol{94}%
}}\left( \{v_{1},v_{2},e_{2},e_{3},e_{2}e_{3}\}\right) $

\strut

$\ \ \ \ \ \ \ \ \ \ \ \ \ \ \ \ \ \ =d\left( \{v_{1},v_{2}\}\right) +\rho
\left( \{e_{2},e_{3},e_{2}e_{3}\}\right) $

\strut

$\ \ \ \ \ \ \ \ \ \ \ \ \ \ \ \ \ \ =0+\left( 1+1+2\right) =\allowbreak 4.$

\strut

Now, consider the element $l=e_{1}e_{2}e_{3}$ in $D_{r}(G^{\symbol{94}}).$
Then

\strut

$\ \ \ \ \ \ \ \ \ \ \ \ \ \ \delta ^{r}(S_{l}^{l})=\{v_{1},l\}$ \ and \ $%
\delta ^{r}(S_{r}^{l})=\{v_{1},l\},$

so,

\ \ \ $\ \ \ \ \ \ \ \ \ \ \ \ \ \ \ \ \ \ \ \delta ^{r}(S_{l}^{l})\cup
\delta ^{r}(S_{r}^{l})=\{v_{1},l\}.$

\strut Therefore,

\strut

$\ \ \ \ \ \ \ \ \ \ \ \ \ \ I_{G_{\Delta }}\left( g_{l}\right) =\mu _{G^{%
\symbol{94}}}\left( \{v_{1},l\}\right) =3.$

\strut \strut \strut \strut \strut
\end{example}

\strut \strut \strut \strut \strut

\strut

\strut

\section{Von Neumann Algebras Induced by Graph Measures}

\strut

\strut

\strut

Throughout this chapter, let $G$ be a finite directed graph and $G^{\symbol{%
94}},$ the shadowed graph of $G$ and let $\Bbb{F}^{+}(G^{\symbol{94}}),$ $%
D(G^{\symbol{94}})$, $\Bbb{G},$ and $D_{r}(G^{\symbol{94}})$ be the free
semigroupoid of $G^{\symbol{94}},$ the diagram set of $G^{\symbol{94}},$ the
graph groupoid of $G$ and the reduced diagram set of $G^{\symbol{94}},$
respectively. Also, for the given algebraic structures, let $\mu _{G^{%
\symbol{94}}},$ $\mu _{\delta },$ $\mu _{\Bbb{G}}$ and $\mu _{r}$ be the
energy measure, the diagram measure, the graph groupoid measure and the
reduced diagram measure of $G,$ respectively.

\strut

\textbf{Notation} \ By $\mu _{G},$ we denote one of the above graph measures
if there is no confusion. $\square $

\strut

\begin{definition}
Let $G$ be the given graph and let $\mu _{G}$ be a graph measure. Define the
graph Hilbert space $H_{G}$ by the space $L^{2}$ $(\mu _{G})$ of all square
integrable $\mu _{G}$-measurable functions.
\end{definition}

\strut \strut

Consider the set $L^{\infty }(\mu _{G})$ of all bounded $\mu _{G}$%
-measurable functions. i.e., $g$ $\in $ $L^{\infty }(\mu _{G})$ if and only
if

\strut

\begin{center}
$\left\| g\right\| _{\infty }\overset{def}{=}\sup \{\left| g(w)\right| :w\in 
\mathcal{G}\}<\infty ,$
\end{center}

\strut

where $\mathcal{G}$ is one of $\Bbb{F}^{+}(G^{\symbol{94}})$, $D(G^{\symbol{%
94}}),$ $\Bbb{G}$ and $D_{r}(G^{\symbol{94}}),$ for the fixed measure $\mu
_{G}.$ Then each element $g$ in $L^{\infty }(\mu _{G})$ can be regarded as a
multiplication operator $M_{g}$ with its symbol $g$ on $H_{G}$ with its
operator norm

\strut

\begin{center}
$\left\| M_{g}\right\| =\left\| g\right\| _{\infty },$ \ for all \ $g$ $\in $
$L^{\infty }(\mu _{G}).$
\end{center}

\strut

It is well-known that the operator algebra $\{M_{g}$ $:$ $g$ $\in $ $%
L^{\infty }(\mu _{G})\},$ as a subalgebra of $B(H_{G}),$ the collection of
all bounded linear operators on $H_{G},$ is a von Neumann algebra.

\strut

\begin{definition}
Let $G$ be a finite directed graph and let $H_{G}$ be the graph Hilbert
space and $B(H_{G}),$ the set of all bounded linear operators on $H_{G}.$
Then the von Neumann algebra $L^{\infty }(\mu _{G})$ $\subset $ $B(H_{G}),$
denoted by $M_{G},$ is called a graph von Neumann algebra.
\end{definition}

\strut \strut \strut

By definition, we can verify the following theorem.

\strut

\begin{theorem}
Let $G_{1}$ and $G_{2}$ be finite directed graphs. Then the graph von
Neumann algebras $M_{G_{1}}$ and $M_{G_{2}}$ are $*$-isomorphic if $G_{1}$
and $G_{2}$ are graph-isomorphic.
\end{theorem}

\strut

\begin{proof}
Assume that $G_{1}$ and $G_{2}$ are graph-isomorphic. Then the graph
measures $\mu _{G_{1}}$ and $\mu _{G_{2}}$ are equivalent. Thus the graph
Hilbert spaces $H_{G_{1}}$ and $H_{G_{2}}$ are isomorphic, since $H_{G_{k}}$ 
$=$ $L^{2}(\mu _{G_{k}}),$ for $k$ $=$ $1,$ $2.$ So, the graph von Neumann
algebras $L^{\infty }(\mu _{G_{1}})$ and $L^{\infty }(\mu _{G_{2}})$ are $*$%
-isomorphic, as $W^{*}$-subalgebras in $B(H_{G_{1}})$ and $B(H_{G_{2}}),$
respectively.
\end{proof}

\strut \strut \strut \strut \strut

\strut \strut

\strut \strut \strut \strut \strut \strut

\strut\textbf{References}

\strut

\begin{quote}
\strut

{\small [1] \ \ A. G. Myasnikov and V. Shapilrain (editors), Group Theory,
Statistics and Cryptography, Contemporary Math, 360, (2003) AMS.}

{\small [2] \ \ D. G. Radcliffe, Rigidity of Graph Products of Groups, Alg
\& Geom. Topology, Vol 3, (2003) 1079 - 1088.}

{\small [3] \ \ D. Voiculescu, Entropy of Random Walks on Groups and the
Macaev Norm, Proc. AMS, Vol 119, 3, (1993), 971 - 977.}

{\small [4] \ \ E. Breuillard and T. Gelander, Cheeger Constant and
Algebraic Entropy of Linear Groups, (2005) Preprint.}

{\small [5] \ \ F. Balacheff, Volum Entropy, Systole and Stable Norm on
Graphs, (2004) Preprint.}

{\small [6] \ \ G. C. Bell, Growth of the Asymptotic Dimension Function for
Groups, (2005) Preprint.}

{\small [7] \ \ I. Cho, Random Variables in a Graph }$W^{*}${\small %
-Probability Space, Ph. D. Thesis, (2005) Univ. of Iowa.}

{\small [8] \ \ I. Cho, Graph Von Neumann Algebras and Graph }$W^{*}${\small %
-Probability Spaces : Crossed Product Approach (2006) Preprint.}

{\small [9] \ \ I. Cho, Operator Theory on Graphs, (2006), In Progress.}

{\small [10] J. Friedman and J-P. Tillich, Calculus on Graphs, (2005)
Preprint.}

{\small [11] J. Stallings, Centerless Groups-An Algebraic Formulation of
Gottlieb's Theorem, Topology, Vol 4, (1965) 129 - 134.}

{\small [12] R. G. Bartle, The Elements of Integration, (1966) 1-st Edition,
John Wiley \& Sons.}

{\small [13] R. Gliman, V. Shpilrain and A. G. Myasnikov (editors),
Computational and Statistical Group Theory, Contemporary Math, 298, (2001)
AMS.}
\end{quote}

\end{document}